\documentclass[a4paper, 11pt]{amsart}

\usepackage{amsmath}
\usepackage{amssymb}
\usepackage{amsthm}
\usepackage{mathtools}
\usepackage{subcaption}
\usepackage[unicode, hidelinks]{hyperref}

\usepackage{tikz}
\usetikzlibrary{calc}

\usepackage[textwidth=16cm, textheight=23cm, centering]{geometry}
\numberwithin{equation}{section}
\allowdisplaybreaks[3]

\newtheorem{theorem}{Theorem}[section]
\newtheorem{proposition}[theorem]{Proposition}
\newtheorem{lemma}[theorem]{Lemma}
\newtheorem{corollary}[theorem]{Corollary}

\theoremstyle{definition}
\newtheorem{definition}[theorem]{Definition}
\newtheorem{example}[theorem]{Example}

\newtheorem{assumption}[theorem]{Assumption}

\theoremstyle{remark}
\newtheorem{remark}[theorem]{Remark}

\newcommand{\tref}[2][]
	{\ifx#1""{Theorem~\textup{\ref{#2}}}%
	\else{Theorems~\textup{\ref{#1}} and~\textup{\ref{#2}}}\fi}
\newcommand{\trefs}[4][]
	{\ifx#1""{Theorems~\textup{\ref{#2}},~\textup{\ref{#3}} and~\textup{\ref{#4}}}%
	\else{Theorems~\textup{\ref{#1}},~\textup{\ref{#2}},~\textup{\ref{#3}} and~\textup{\ref{#4}}}\fi}
\newcommand{\pref}[2][]
	{\ifx#1""{Proposition~\textup{\ref{#2}}}%
	\else{Propositions~\textup{\ref{#1}} and~\textup{\ref{#2}}}\fi}
\newcommand{\prefs}[4][]
	{\ifx#1""{Propositions~\textup{\ref{#2}},~\textup{\ref{#3}} and~\textup{\ref{#4}}}%
	\else{Propositions~\textup{\ref{#1}},~\textup{\ref{#2}},~\textup{\ref{#3}} and~\textup{\ref{#4}}}\fi}
\newcommand{\lref}[2][]
	{\ifx#1""{Lemma~\textup{\ref{#2}}}%
	\else{Lemmas~\textup{\ref{#1}} and~\textup{\ref{#2}}}\fi}
\newcommand{\lrefs}[4][]
	{\ifx#1""{Lemma~\textup{\ref{#2}}, \textup{\ref{#3}} and~\textup{\ref{#4}}}%
	\else{Lemmas~\textup{\ref{#1}}, \textup{\ref{#2}}, \textup{\ref{#3}} and \textup{\ref{#4}}}\fi}
\newcommand{\corref}[2][]
	{\ifx#1""{Corollary~\textup{\ref{#2}}}%
	\else{Corollaries~\textup{\ref{#1}} and~\textup{\ref{#2}}}\fi}
\newcommand{\exref}[2][]
	{\ifx#1""{Example~\textup{\ref{#2}}}%
	\else{Examples~\textup{\ref{#1}} and~\textup{\ref{#2}}}\fi}
\newcommand{\rref}[2][]
	{\ifx#1""{Remark~\textup{\ref{#2}}}%
	\else{Remark~\textup{\ref{#1}} and~\textup{\ref{#2}}}\fi}
\newcommand{\aref}[2][]
	{\ifx#1""{Assumption~\textup{\ref{#2}}}%
	\else{Assumption~\textup{\ref{#1}} and~\textup{\ref{#2}}}\fi}
\newcommand{\secref}[2][]
	{\ifx#1""{Section~\textup{\ref{#2}}}%
	\else{Sections~\textup{\ref{#1}} and~\textup{\ref{#2}}}\fi}
\newcommand{\figref}[2][]
	{\ifx#1""{Figure~\textup{\ref{#2}}}%
	\else{Figure~\textup{\ref{#1}} and~\textup{\ref{#2}}}\fi}
\newcommand{\itemref}[1]{\textup{(\ref{#1})}}

\newcommand{\RealNum}{\mathbf{R}}
\newcommand{\NaturalNum}{\mathbf{N}}
\newcommand{\Integers}{\mathbf{Z}}

\newcommand{\even}{\mathfrak{e}}
\newcommand{\odd}{\mathfrak{o}}
\newcommand{\op}{\mathrm{op}}
\newcommand{\onevector}{\mathbf{1}}
\newcommand{\zerovector}{\mathbf{0}}
\DeclareMathOperator{\trace}{tr}
\newcommand{\diag}{\mathrm{diag}}

\newcommand{\sigmaField}{\mathcal{F}}
\newcommand{\prob}{\boldsymbol{P}}
\newcommand{\expect}{\boldsymbol{E}}

\newcommand{\covariance}{\textup{Cov}}

\newcommand{\indicator}[1]{\mathtt{1}_{#1}}
\newcommand{\setindicator}[1]{\mathtt{1}_{\{#1\}}}

\begin{document}

\title[ERWs on Coverings of Dipole Graphs]{Elephant Random Walks on Coverings of Dipole Graphs}

\author[N. Naganuma]{Nobuaki Naganuma}
\address[N. N.]{Faculty of Advanced Science and Technology, Kumamoto University.
2-39-1 Kurokami, Chuo-ku, Kumamoto, 860-8555, JAPAN}
\email{\url{naganuma@kumamoto-u.ac.jp}}

\author[K. Yura]{Kaito Yura}
\address[K. Y.]{Graduate School of Science and Technology, Kumamoto University.
2-39-1 Kurokami, Chuo-ku, Kumamoto, 860-8555, JAPAN}
\email{\url{kaity_yura@kumadai.jp}}

\thanks{
The authors thank Professors Masato Takei, Syota Esaki, and Ryuya Namba for fruitful discussions
and Professor Hiroshi Kawabi for bringing reference \cite{Ureshino2021} to their attention.
They are also grateful to the anonymous referee for valuable comments.
The first named author is partially supported by JSPS KAKENHI Grant Number 22K13932.
During the preparation of this manuscript, the authors used ChatGPT (OpenAI). All mathematical content was reviewed and verified by the authors.}

\begin{abstract}
	In the present paper, we introduce and analyze elephant random walks (ERWs)
	on bipartite periodic lattices arising as coverings of dipole graphs.
	We focus on lattices whose admissible step directions 
	in the two parts of the bipartition are negatives of each other and disjoint. 
	On such graphs, we define an ERW in which each step 
	is chosen by referring to the entire history of the walk.
	The ERW on the hexagonal lattice is a prototypical example of our model.
	The definition and asymptotic analysis of such ERWs are not straightforward
	because both depend strongly on the underlying geometric structure.
	
	Our analysis is based on a combination of the P\'olya-type urn techniques and the martingale approach, 
	two standard methods for analyzing ERWs.
	We find that the counting process of the ERW forms a P\'olya-type urn with two-periodic generating matrices.
	By analyzing such urn models,
	we show the strong law of large numbers for the counting process.
	Combining the result for the counting process with the martingale approach,
	we derive the strong laws of large numbers with non-standard normalizations
	and central limit theorems
	for the position process of the ERW in the diffusive and critical regimes,
	as well as almost sure and $L^2$ scaling limits in the superdiffusive regime.
\end{abstract}

\maketitle

\setcounter{tocdepth}{1}
\tableofcontents

\section{Introduction}
\subsection{Background}
Random walks are fundamental stochastic processes in probability theory and have been extensively studied in various contexts.
Starting from the simple random walk on the one-dimensional integer lattice,
it is natural to consider its higher-dimensional extension on $\Integers^d$ and on more general graphs.
In such settings, the geometric structure of the underlying space influences the long-term behavior of the walk,
leading to rich phenomena and a wide range of asymptotic behaviors.

Beyond simple random walks, models incorporating memory effects have attracted significant attention,
especially in statistical physics and related areas.
Sch\"utz and Trimper \cite{SchutzTrimper2004ERW} introduced the elephant random walk (ERW),
a discrete-time stochastic process with long-term memory on $\Integers$,
to investigate how dependence on the past affects scaling behavior and induces a phase transition
from diffusive to superdiffusive regimes.
A key feature of the ERW is a memory parameter that separates diffusive, critical, and superdiffusive regimes.

The ERW can often be analyzed through the P\'olya-type urn techniques and the martingale approach.
Baur and Bertoin \cite{BaurBertoin2016ERWUrn} established a connection between the ERW and P\'olya-type urns,
enabling the use of urn theory to analyze functional limit theorems for the ERW.
Coletti, Gava and Sch\"utz \cite{ColettiGavaSchutz2017CLT,ColettiGavaSchutz2017SIP}
applied the martingale approach to the ERW
and derived various laws of large numbers, central limit theorems and a strong invariance principle.
Bercu \cite{Bercu2018MartingaleERW} also used the martingale approach and derived the law of iterated logarithm for the ERW.

The ERW has been studied in various settings.
One natural extension is to consider ERWs on higher-dimensional lattices
(e.g., 
\cite{
	Bercu2025MERWWithStops,
	BercuLaulin2019MERW, 
	Bertenghi2022FunctionalMERW,
	ChenLaulin2023MARW,
	CurienLaulin2024PlaneRecurrenceMERW,
	GhoshDhillonKataria2026MovesMERWS,
	GonzalezNavarrete2020RandomTendencyMERW,
	GuerinLaulinRaschel2023FixedPointSuperdiffusiveMERW,
	Marquioni2019CoupledMemoryMERW, 
	Qin2025RecTransMERW}).
Bercu and Laulin \cite{BercuLaulin2019MERW} initiated the study of the multi-dimensional ERW (MERW)
and established limit theorems for the MERW via the martingale approach.
Bertenghi \cite{Bertenghi2022FunctionalMERW} studied functional limit theorems for the MERW
by P\'olya-type urn techniques.
More recently, Bercu \cite{Bercu2025MERWWithStops} introduced the MERW with stops and analyzed its asymptotic behavior.
Apart from extensions to higher dimensions, Mukherjee \cite{Mukherjee2025CayleyTree}, Shibata \cite{Shibata2025ArXivPeriodic}
and Ureshino \cite{Ureshino2021}
studied ERWs on graphs.

\subsection{A prototypical example}\label{sec940290329}
The ERW on the hexagonal lattice is a prototypical example to which our main results apply.
We explain here why the hexagonal lattice is a natural and important example
from the viewpoint of determining the next step from the history of the walk,
which is the basic idea of ERWs.

\begin{figure}[t]
	\begin{minipage}[b]{0.32\linewidth}
	    \centering
		\begin{tikzpicture}[scale=0.8]
			\clip (-2.8,-2.8) rectangle (2.8,2.8);
			
			\coordinate (vectorU) at ($(1,0)$); 
			\coordinate (vectorV) at ($(0,1)$); 
			\coordinate (vectorW) at ($(-1,0)$);

			\foreach \i in {-5,...,5}{
				\foreach \j in {-5,...,5}{
					\coordinate (A) at ($\i*(vectorU)-\i*(vectorW)+\j*(vectorV)-\j*(vectorW)$);
					\draw (A) -- ($(A)+(vectorU)$);
					\draw (A) -- ($(A)+(vectorV)$);
					\draw (A) -- ($(A)+(vectorW)$);
					\draw (A) -- ($(A)+(vectorU)-(vectorV)+(vectorW)$);
				}
			}

			\foreach \i in {-5,...,5}{
				\foreach \j in {-5,...,5}{
					\coordinate (A) at ($\i*(vectorU)-\i*(vectorW)+\j*(vectorV)-\j*(vectorW)$);
					\fill (A) circle (3pt);
					\filldraw[fill=white] ($(A)+(vectorU)$) circle (3pt);
				}
			}
		\end{tikzpicture}
		\subcaption{Square}
		\label{figSquare}
	\end{minipage}
	\begin{minipage}[b]{0.32\linewidth}
	    \centering
		\begin{tikzpicture}[scale=0.8]
			\clip (-2.3,-2.8) rectangle (3.3,2.8);
			
			\coordinate (vectorU) at ($(1,0)$); 
			\coordinate (vectorV) at ($(-1/2,1.732050/2)$); 
			\coordinate (vectorW) at ($(-1/2,-1.732050/2)$);

			\foreach \i in {-5,...,5}{
				\foreach \j in {-5,...,5}{
					\coordinate (A) at ($\i*(vectorU)-\i*(vectorV)+\j*(vectorU)-\j*(vectorW)$);
					\draw (A) -- ($(A)+(vectorU)$);
					\draw (A) -- ($(A)+(vectorV)$);
					\draw (A) -- ($(A)+(vectorW)$);
				}
			}

			\foreach \i in {-5,...,5}{
				\foreach \j in {-5,...,5}{
					\coordinate (A) at ($\i*(vectorU)-\i*(vectorV)+\j*(vectorU)-\j*(vectorW)$);
					\fill (A) circle (3pt);
					\filldraw[fill=white] ($(A)+(vectorU)$) circle (3pt);
				}
			}
		\end{tikzpicture}
		\subcaption{Hexagonal}
		\label{figHexagonal}
	\end{minipage}
	\begin{minipage}[b]{0.32\linewidth}
		\centering
		\begin{tikzpicture}[scale=0.8]
			\clip (-2.8,-2.8) rectangle (2.8,2.8);
			
			\coordinate (vectorU) at ($(1,0)$); 
			\coordinate (vectorV) at ($(0,1)$); 
			\coordinate (vectorW) at ($(-1,0)$);

			\foreach \i in {-5,...,5}{
				\foreach \j in {-5,...,5}{
					\coordinate (A) at ($\i*(vectorU)-\i*(vectorW)+\j*(vectorV)-\j*(vectorW)$);
					\draw (A) -- ($(A)+(vectorU)$);
					\draw (A) -- ($(A)+(vectorV)$);
					\draw (A) -- ($(A)+(vectorW)$);
				}
			}

			\foreach \i in {-5,...,5}{
				\foreach \j in {-5,...,5}{
					\coordinate (A) at ($\i*(vectorU)-\i*(vectorW)+\j*(vectorV)-\j*(vectorW)$);
					\fill (A) circle (3pt);
					\filldraw[fill=white] ($(A)+(vectorU)$) circle (3pt);
				}
			}
		\end{tikzpicture}
		\subcaption{Brick-wall}
		\label{figBrickWall}
	\end{minipage}
	\caption{Some periodic planar lattices.}
	\label{fig:lattices-overview}
\end{figure}
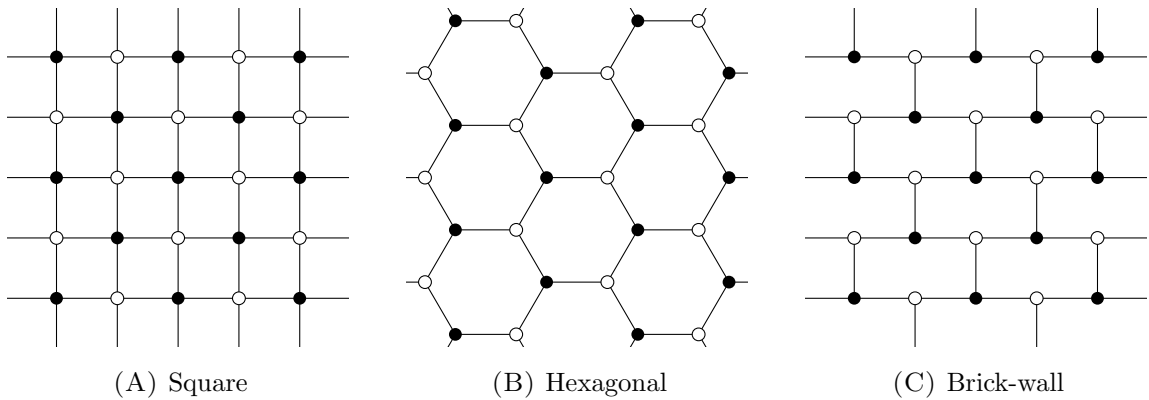

Typical examples of planar lattices include the square, hexagonal, and brick-wall lattices.
For these lattices, we focus on the local structure at each vertex.
More precisely, for a vertex $x$, let $V_x$ denote the set of directions of edges leaving $x$.
For the examples above, the sets $V_x$ are as follows.
\begin{itemize}
	\item	The square lattice: $V_x  = \{\pm\mathbf{e}_1,\pm\mathbf{e}_2\}$ for all vertices $x$ in \figref{figSquare}.
	\item	The hexagonal lattice:
			\begin{gather*}
				\begin{aligned}
					V_x
					&=
						\left\{
							+\mathbf{e}_1,
							\frac{1}{2}
							\big(
								-\mathbf{e}_1\pm \sqrt{3}\mathbf{e}_2
							\big)
						\right\}
					&
					&
					\text{for $x$ at the black dots in \figref{figHexagonal},}\\
					V_x
					&=
						\left\{
							-\mathbf{e}_1,
							\frac{1}{2}
							\big(
								+\mathbf{e}_1\mp \sqrt{3}\mathbf{e}_2
							\big)
						\right\}
					&
					&
					\text{for $x$ at the white dots in \figref{figHexagonal}.}
				\end{aligned}
			\end{gather*}

	\item	The brick-wall lattice:
			\begin{gather*}
				\begin{aligned}
					V_x &= \{\pm\mathbf{e}_1,+\mathbf{e}_2\}
					&
					&
					\text{for $x$ at the black dots in \figref{figBrickWall},}\\
					V_x &= \{\mp\mathbf{e}_1,-\mathbf{e}_2\}
					&
					&
					\text{for $x$ at the white dots in \figref{figBrickWall}.}
				\end{aligned}
			\end{gather*}
\end{itemize}
Here, $\mathbf{e}_1,\mathbf{e}_2$ denote the standard basis of $\RealNum^2$.
In \figref{fig:lattices-overview} and subsequent figures, 
a coordinate system with one of the black dots $\bullet$ as the origin is assumed;
however, the axes are omitted for simplicity.
Here and throughout, after fixing an appropriate unit length,
we express vectors in terms of their coordinates.
These examples lead to three observations.
First, there are at most two distinct types of local structures, denoted by $V_\bullet$ and $V_\circ$.
Indeed, the square lattice consists of a single type of local structure, with $V_\bullet=V_\circ$,
whereas the hexagonal and brick-wall lattices consist of two types of local structures, with $V_\bullet\neq V_\circ$.
Second, $V_\bullet= -V_\circ$,
which follows from the bipartite structure of the lattices.
Third, we obtain the following trichotomy:
the identical case $V_\bullet = V_\circ$,
the disjoint case $V_\bullet \cap V_\circ = \emptyset$,
and the partially overlapping case $V_\bullet \cap V_\circ \neq \emptyset$ and $V_\bullet \neq V_\circ$.

We consider the three observations from the viewpoint of determining the next step from the history of the walk.
The first observation leads to a difficulty in defining ERWs on graphs.
In the case where every vertex has the same local structure, that is, $V_\bullet = V_\circ$,
ERWs can be defined and analyzed in the same manner as MERWs.
On the other hand, in the case where there are two types of local structures, that is, $V_\bullet \neq V_\circ$,
the definition is no longer straightforward.
More concretely, when the elephant is at a vertex with $V_\bullet$ (resp. $V_\circ$)
and remembers a direction in $V_\circ$ (resp. $V_\bullet$),
it is not clear how to define the corresponding move.
However, the second observation $V_\bullet= -V_\circ$ is a key ingredient in overcoming this difficulty.
The condition gives a natural correspondence between elements of $V_\bullet$ and those of $V_\circ$
and ensures that we can define the transition probabilities of the ERW.
See \secref{sec9092423242092} for the definition.
The third observation, together with the discussion below,
leads us to focus on the disjoint case $V_\bullet \cap V_\circ = \emptyset$.
As stated above, the identical case $V_\bullet = V_\circ$ can be treated
similarly to Bercu and Laulin \cite{BercuLaulin2019MERW}.
We therefore turn to the non-identical case.
Among them, the disjoint case is the simplest. 
Indeed, the parity of time alone determines to which local structure a step belongs. 
In the partially overlapping case, by contrast, this is no longer sufficient:
we must additionally keep track of how many times steps in the overlap have been taken, 
which makes the situation more complicated.
In \secref{secConcludingRemarks}, we give more information 
about the identical and partially overlapping cases.

\subsection{Overview of main results}
Using $V_\bullet= -V_\circ$, we can define the ERW on the hexagonal lattice
so that the elephant refers to the entire history of the walk,
even though the admissible directions vary with its location. 
This is a key difference from the models of Ureshino \cite{Ureshino2021}, Shibata \cite{Shibata2025ArXivPeriodic}, 
and Mukherjee \cite{Mukherjee2025CayleyTree}.
Ureshino and Shibata independently studied ERWs on periodic lattices,
whose transition mechanism is defined separately for each local structure.
More precisely, when the elephant is currently at a vertex with $V_\bullet$ (resp.\ $V_\circ$),
it refers only to the past steps taken at vertices with $V_\bullet$ (resp.\ $V_\circ$).
Furthermore, Mukherjee introduces an ERW on infinite Cayley trees,
where the elephant is defined in graph-theoretic terms (without an embedding into $\RealNum^d$).
In his model, the labelled local neighborhood is identical at every vertex.
Thus, in both models, the next step is determined by
past steps taken at vertices with the same local structure as the current vertex.
In this sense, these models can be viewed as an extension of MERWs.
In contrast, we define the ERW so that the elephant refers
to the entire history of the walk, 
including steps taken at the other local structure.

We next explain our setting and method of analysis.
To place this example in a broader context, we proceed through three levels of abstraction.
The first is the ERW on the hexagonal lattice, 
which serves as our primary motivating example.
We then introduce ERWs on periodic realizations of 
coverings of dipole graphs under suitable assumptions on the edge directions.
Finally, we formulate the model under abstract assumptions
on the state space and the step directions; this abstract framework is presented in \secref{sec:setting}.
The key structural feature that emerges from the hexagonal lattice example is the relation $V_\bullet=-V_\circ$.
In fact, this is the only geometric property used in the subsequent analysis.

We then describe how we analyze the ERW in this general setting.
First, we construct a counting process that counts the number of times each step
has been taken up to time $n$.
This counting process forms a P\'olya-type urn with two-periodic generating matrices.
Inspired by the results for such urns established by Yan-Cheng-Bai \cite{YanChengBai2006},
we derive the standard strong law of large numbers for the counting and position processes.
Next, we characterize the phase transition
by analyzing the centered second moment of the position process.
As a result, we find that the ERW exhibits three regimes: diffusive, critical, and superdiffusive.
In each regime, we derive the corresponding limit theorems by combining the martingale approach
and the limit theorems for the associated urn model.
In diffusive and critical regimes, we obtain 
the strong laws of large numbers with non-standard normalizations and central limit theorems.
In the superdiffusive regime, we prove the scaling limit of the position process.

\subsection{Organization and notation}
The remainder of the paper is organized as follows.
In \secref{sec:setting}, we review coverings of dipole graphs
and then formulate the ERW in the abstract setting used throughout the paper.
We also state the main limit theorems
for the counting and position processes of the ERW.
In \secref{sec320492094209410}, we analyze the counting process
and show the corresponding limit theorems. 
The standard strong law of large numbers for the position process is also proved.
\secref{sec58934850298092} is devoted to calculations of
the conditional expectation of one-step increments,
which play an important role in the martingale approach.
In \secref{sec:Centered-second-moment}, we analyze the centered second moment
to characterize the phase transition.
In \secref{sec:limit-theorems}, we construct a suitable martingale associated with the ERW 
and prove the non-standard strong law of large numbers 
and central limit theorems in the diffusive and critical regimes.
Furthermore, we establish almost sure and $L^2$ scaling limits
in the superdiffusive regime.
\secref{secConcludingRemarks} is devoted to concluding remarks,
in which we comment on the identical and partially overlapping cases.
In Appendix~\ref{appendix}, we collect limit theorems on martingales.

We end the introduction with notation used throughout the paper.
\begin{itemize}
	\item	Let $\{\mathbf{e}_i\}_{i=1}^N$ be the standard basis in $\RealNum^N$.
	\item 	Write $\onevector_N = (1,\dots,1)^\top \in \RealNum^N$ and $\zerovector_N = (0,\dots,0)^\top \in \RealNum^N$.
	\item	For matrices $A$, we define the operator norm by $\|A\|_\op = \sup_{v\in\RealNum^N,|v|=1}|Av|$.
	\item   Let $\indicator{\odd}(n)$ and $\indicator{\even}(n)$ be the indicator functions of parity of $n\in\Integers$ defined by
			$\indicator{\odd}(n) = 1$ if $n$ is odd, and $0$ otherwise, and
			$\indicator{\even}(n) = 1$ if $n$ is even, and $0$ otherwise.
	\item  	Let $(a_n)$ and $(b_n)$ be real sequences with $b_n > 0$.  
			We write
			$a_n = O(b_n)$ if there exists $C > 0$ and $N \in \NaturalNum$ such that $|a_n| \le C b_n$ for all $n \ge N$.
			$a_n = o(b_n)$ if $\lim_{n \to \infty} a_n/b_n = 0$.
			$a_n \sim b_n$ if $\lim_{n \to \infty} a_n/b_n = 1$.
	\item	Denote by $\Gamma(\cdot)$ the Gamma function,
			which is defined by
			$
				\Gamma(x)
				=
					\int_0^\infty
						e^{-t}t^{x-1}\,
						dt
			$, $x>0$.
			Note that, for $a,b \in \RealNum$, by Stirling's formula we have
			\begin{align}\label{eq:Gamma-asymptotic}
				\frac{\Gamma(x + a)}{\Gamma(x + b)}
				\sim
					x^{a-b}.
			\end{align}
	\item	Denote by $\mathcal{N}_d(\zerovector_d,\Sigma)$ 
			the $d$-dimensional Gaussian distribution with mean $\zerovector_d$ and covariance matrix $\Sigma$.
	\item	Denote by $\sharp A$ the cardinality of a set $A$.
	\item   Let $C$ be a generic positive constant whose value may change from line to line.
\end{itemize}

\section{Setting of elephant random walks}\label{sec:setting}
In this section, we define the ERW and present the main results of this paper.
In \secref{sec487239842093841}, we explain the observations concerning $V_\bullet$ and $V_\circ$ 
in the framework of topological crystallography 
(Kotani-Sunada \cite{KotaniSunada2000Jacobian, KotaniSunada2001Standard,KotaniSunada2003} and Sunada \cite{Sunada2012, Sunada2013Book}).
More precisely, we review graphs, topological crystals, and periodic realizations. 
In particular, we explain coverings of dipole graphs and give several examples.
We introduce them only to help the reader understand the overall picture
and concrete examples; the subsequent analysis does not use them.
In \secref{sec9092423242092}, we define the ERW under abstract assumptions 
that include the graphs introduced in \secref{sec487239842093841} as examples.
In \secref{sec5409480923184092}, we state the main results.
The main theorems are proved under the abstract assumptions.
Readers who wish to start with the ERW on the hexagonal lattice 
may first read \secref[sec9092423242092]{sec5409480923184092}.

\subsection{Periodic realizations and dipole graphs}\label{sec487239842093841}
We set up some notation related to graphs, topological crystals and periodic realizations 
based on Sunada \cite{Sunada2012, Sunada2013Book}.
We recall some definitions from them.
\begin{definition}
	A (non-oriented) graph is an ordered pair $\mathcal{G}=(\mathcal{X},\mathcal{E})$ 
	of disjoint sets, together with maps $o,t:\mathcal{E}\to \mathcal{X}$ 
	and an involution $e\mapsto e^{-1}$ on $\mathcal{E}$
	such that
	\begin{align*}
		(e^{-1})^{-1}=e,\qquad
		o(e^{-1})=t(e),\qquad
		t(e^{-1})=o(e)
		\qquad (e\in\mathcal{E}).
	\end{align*}
	Elements of $\mathcal{X}$ are called vertices
	and elements of $\mathcal{E}$ are called directed edges.
  	For each edge $e\in\mathcal{E}$,
	the vertices $o(e)$ and $t(e)$ are called the origin and terminus of $e$, respectively.
	
	The graph $\mathcal{G}$ is called a finite graph if $\mathcal{X}$ and $\mathcal{E}$ are both finite sets.
	The graph $\mathcal{G}$ is said to be locally finite if $\sharp \mathcal{E}_x<\infty$ for all $x\in \mathcal{X}$.
	Here, we write $\mathcal{E}_x=\{e\in\mathcal{E}\,|\,o(e)=x\}$ for $x\in \mathcal{X}$.
\end{definition}

Next, we introduce the notion of topological crystals.
\begin{definition}
	Let $\mathcal{G}=(\mathcal{X},\mathcal{E})$ and 
	$\mathcal{G}_0=(\mathcal{X}_0,\mathcal{E}_0)$ be connected, locally finite graphs.
	\begin{enumerate}
		\item	A graph morphism $\omega : \mathcal{G} \to \mathcal{G}_0$ is called a covering map
				if the following holds.
				\begin{enumerate}
					\item	$\omega : \mathcal{X} \to \mathcal{X}_0$
							is surjective
					\item	for every vertex $x \in \mathcal{X}$, 
							$\omega|_{\mathcal{E}_x} : \mathcal{E}_x \to (\mathcal{E}_0)_{\omega(x)}$
							is bijective.
				\end{enumerate}
				Then, $\mathcal{G}_0$ is called the base graph and
				the group $\Gamma=\{g\in \mathop{\mathrm{Aut}}(\mathcal{G})\,|\,\omega\circ g=\omega\}$
				is called the covering transformation group of the covering map $\omega$.
		\item	A covering map $\omega : \mathcal{G} \to \mathcal{G}_0$ is said to be regular
				if, for any $x,y\in\mathcal{X}$ with $\omega(x)=\omega(y)$,
				there exists an element $g\in\Gamma$ such that $y=gx$.
				Here, we write $gx=g(x)$ for notational simplicity.
	\end{enumerate}
\end{definition}

\begin{definition}
	Let $\mathcal{G}$ be a connected, locally finite graph
	and $\mathcal{G}_0$ a connected finite graph.
	We say that $\mathcal{G}$ is a $d$-dimensional topological crystal over $\mathcal{G}_0$ 
	if there exists a regular covering map $\omega : \mathcal{G} \to \mathcal{G}_0$
	whose covering transformation group $\Gamma$ is a free abelian group of rank $d$.
\end{definition}

Next we consider a periodic realization of $\mathcal{G}=(\mathcal{X}, \mathcal{E})$ in $\RealNum^d$.
\begin{definition}
	Let $\mathcal{G}= (\mathcal{X}, \mathcal{E})$ be a $d$-dimensional topological crystal 
	with a covering transformation group $\Gamma$,
	and let $\mathcal{G}_0=(\mathcal{X}_0,\mathcal{E}_0)$ be its base graph.
	A piecewise linear map $\Phi:\mathcal{G}\to\RealNum^d$ is called a periodic realization
	if there exists an injective homomorphism $\rho:\Gamma\to\RealNum^d$ 
	whose image is a lattice in $\RealNum^d$ such that
	\begin{align*}
		\Phi(g x)
		=
		\Phi(x)+\rho(g)
		\qquad (x\in \mathcal{X},\ g\in \Gamma).
	\end{align*}
\end{definition}
We note that we can define $\Phi:\mathcal{E}\to\RealNum^d\times\RealNum^d$ by
$
	\Phi(e)
	=
		(\Phi(o(e)),\Phi(t(e)))
$.
Then, since $o(ge)=go(e)$, we have
\begin{align*}
	\Phi(ge)
	=
		(\Phi(o(ge)),\Phi(t(ge)))
	=
		(\Phi(o(e))+\rho(g),\Phi(t(e))+\rho(g))
	=
		\Phi(e)+\rho(g).
\end{align*}
Here, given $a=(a_1,a_2)\in\RealNum^d\times\RealNum^d$ and $\mathbf{v}\in\RealNum^d$,
we wrote $a+\mathbf{v}=(a_1+\mathbf{v},a_2+\mathbf{v})$.

\begin{figure}
    \begin{tikzpicture}
		\draw (0,0) -- (5,0);
		\draw (0,0) to[bend left=80]  (5,0);
		\draw (0,0) to[bend right=80] (5,0);

		\fill (0,0) circle (3pt);
		\filldraw[fill=white] (5,0) circle (3pt);
	\end{tikzpicture}
	\caption{Dipole graph with three edges}
	\label{fig:3-edges-dipole-graph}
\end{figure}
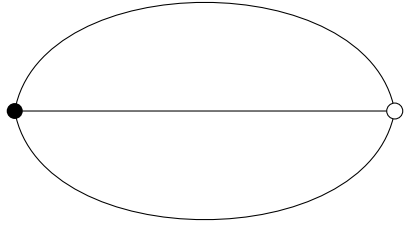

Hereafter, we consider a $d$-dimensional 
topological crystal over dipole graphs and realize it in $\RealNum^d$.
Examples of such realizations in $\RealNum^2$
include the square, hexagonal and brick-wall lattices.
We begin by defining dipole graphs. See \figref{fig:3-edges-dipole-graph}.
Let $m\in\NaturalNum$.
The dipole graph $\mathcal{D}_m = (\mathcal{X}_m, \mathcal{E}_m)$
is defined to be the finite graph with vertex set 
$\mathcal{X}_m = \{\bullet,\circ\}$ and edge set
$\mathcal{E}_m = \{e_1, \hat{e}_1,\dots,e_m, \hat{e}_m\}$.
Maps $o$ and $t$ are defined by
\begin{align*}
    o(e_i) = \bullet, \quad t(e_i) = \circ,
    \quad
    o(\hat{e}_i) = \circ, \quad t(\hat{e}_i) = \bullet
    \qquad (i=1,2,\dots,m).
\end{align*}
Then the involution on $\mathcal{E}_m$ is given by
\begin{align*}
    e_i^{-1} = \hat{e}_i,
    \quad
    \hat{e}_i^{-1} = e_i
    \qquad (i=1,2,\dots,m).
\end{align*}
Thus we may write 
$\mathcal{E}_m = \{e_1, e_1^{-1},\dots,e_m,e_m^{-1}\}$,
where the involution satisfies
\begin{align*}
    (e_i^{-1})^{-1} = e_i,
    \quad
    o(e_i^{-1}) = t(e_i),
    \quad
    t(e_i^{-1}) = o(e_i)
    \qquad (i=1,2,\dots,m).
\end{align*}

Let $\mathcal{G}=(\mathcal{X},\mathcal{E})$ be
a $d$-dimensional topological crystal over a dipole graph $\mathcal{D}_m$.
Its covering map and covering transformation group are denoted by $\omega$ and $\Gamma$, respectively.
Furthermore, $\Phi:\mathcal{G}\to\RealNum^d$ stands for the realization map.
For simplicity, we assume $\Phi$ is injective, which implies that all vertices are mapped to distinct points.
Under this notation, our purpose is to study the graph $\Phi(\mathcal{G})=(\Phi(\mathcal{X}),\Phi(\mathcal{E}))$,
where
\begin{align*}
	\Phi(\mathcal{X})
	&=
		\{\Phi(x')\in\RealNum^d\,|\,x'\in \mathcal{X}\},
	&
	\Phi(\mathcal{E})
	&=
		\{\Phi(e')\in\RealNum^d\times \RealNum^d\,|\,e'\in \mathcal{E}\}.
\end{align*}
The origin and terminus of $e\in \Phi(\mathcal{E})$ are denoted
by $o(e)$ and $t(e)$, respectively.
By a slight abuse of notation, we use the same symbol for simplicity.
To study $\Phi(\mathcal{G})$, we introduce some additional notation.
Set
\begin{align*}
	[\bullet]_{\mathcal{X}}
	&=
		\omega^{-1}(\{\bullet\})
	=
		\{x'\in \mathcal{X}\,|\,\omega(x')=\bullet\},
	&
	[\circ]_{\mathcal{X}}
	&=
		\omega^{-1}(\{\circ\})
	=
		\{x'\in \mathcal{X}\,|\,\omega(x')=\circ\},\\
	[\bullet]
	&=
		\Phi([\bullet]_{\mathcal{X}})
	=
		\{\Phi(x')\in \RealNum^d\,|\,x'\in[\bullet]_\mathcal{X}\},
	&
	[\circ]
	&=
		\Phi([\circ]_{\mathcal{X}})
	=
		\{\Phi(x')\in \RealNum^d\,|\,x'\in[\circ]_\mathcal{X}\}.
\end{align*}
We set
$
	V_x
	=
		\{ \mathbf{v}(e)\in\RealNum^d\,|\,e\in \Phi(\mathcal{E})_x\}
$
for $x\in \Phi(\mathcal{X})$,
where
$
	\mathbf{v}(e)
	=
		t(e)-o(e)
	\in
		\RealNum^d
$
for $e\in\Phi(\mathcal{E})$
and 
$\Phi(\mathcal{E})_x=\{e\in \Phi(\mathcal{E})\,|\,o(e)=x\}$ for $x\in \Phi(\mathcal{X})$.

From the definition, we have $\mathcal{X}=[\bullet]_{\mathcal{X}}\cup [\circ]_{\mathcal{X}}$
and $[\bullet]_{\mathcal{X}}\cap [\circ]_{\mathcal{X}}=\emptyset$.
For all $e'\in\mathcal{E}$, we see 
$o(e')\in [\bullet]_{\mathcal{X}}$ (resp.\,$[\circ]_{\mathcal{X}}$)
if and only if $t(e')\in [\circ]_{\mathcal{X}}$ (resp.\,$[\bullet]_{\mathcal{X}}$).
From the definition, we have
$
	[\bullet]\cup[\circ]
	=
		\{\Phi(x')\in \RealNum^d\,|\,x'\in [\bullet]_{\mathcal{X}}\cup [\circ]_{\mathcal{X}}\}
	=
		\{\Phi(x')\in \RealNum^d\,|\,x'\in \mathcal{X}\}
	=
		\Phi(\mathcal{X})
$.
Since $\Phi$ is injective, we have $[\bullet]\cap[\circ]=\emptyset$.
We also see that
$o(e)\in [\bullet]$ (resp.\,$[\circ]$)
if and only if $t(e)\in [\circ]$ (resp.\,$[\bullet]$) for all $e\in\Phi(\mathcal{E})$.
Furthermore, we obtain the next proposition.
\begin{proposition}\label{prop4892384092}
	The following holds.
	\begin{enumerate}
		\item	We have $V_x=V_y$ for all $x,y\in[\bullet]$,
				and, hence, we can write $V_\bullet=V_x$ for $x\in[\bullet]$.
				The same assertion holds with $\bullet$ replaced by $\circ$.
		\item	$x+\mathbf{v}\in [\circ]$ for all $x\in[\bullet]$ and $\mathbf{v}\in V_\bullet$.
				The same assertion holds with $\bullet$ and $\circ$ interchanged.
		\item	$V_\bullet=-V_\circ$.
	\end{enumerate}
\end{proposition}
\begin{proof}
	(1)
	Let $x,y\in [\bullet]$.
	Then there exist unique $x',y'\in [\bullet]_\mathcal{X}$
	such that $x=\Phi(x')$ and $y=\Phi(y')$  due to the injectivity of $\Phi$.
	In addition, the regularity of $\omega$ implies there exists $g\in\Gamma$ such that $y'=gx'$.
	We have
	\begin{align*}
		\Phi(\mathcal{E})_y
		=
			\Phi(\mathcal{E}_{y'})
		=
			\Phi(\mathcal{E}_{gx'})
		=
			\Phi(g\mathcal{E}_{x'})
		=
			\Phi(\mathcal{E}_{x'})
			+
			\rho(g)
		=
			\Phi(\mathcal{E})_x
			+
			\rho(g),
	\end{align*}
	which implies the assertion.
	Here, given $A\subset\RealNum^d\times \RealNum^d$ and $\mathbf{v}\in\RealNum^d$,
	we wrote $A+\mathbf{v}=\{a+\mathbf{v}\in \RealNum^d\times \RealNum^d\,|\,a\in A\}$.
	We see the assertion for $x,y\in [\circ]$
	in the same way.

	(2)
	Let $x\in[\bullet]$ and $\mathbf{v}\in V_\bullet$.
	Then, since $\mathbf{v}\in V_\bullet=V_x$,
	there exists $e\in\Phi(\mathcal{E})_x$ such that $\mathbf{v}=\mathbf{v}(e)$.
	Noting the unique existence $x'\in [\bullet]_\mathcal{X}$ such that $x=\Phi(x')$ due to the injectivity of $\Phi$
	and $e'\in\mathcal{E}_{x'}$ such that  $e=\Phi(e')$, we have
	\begin{align*}
		\mathbf{v}
		=
			\mathbf{v}(e)
		=
			t(e)-o(e)
		=
			t(\Phi(e'))-o(\Phi(e'))
		=
			\Phi(t(e'))-\Phi(o(e')),
	\end{align*}
	which implies
	$
		x+\mathbf{v}=\Phi(t(e'))
	$.
	Since $t(e')\in [\circ]_{\mathcal{X}}$, we see $x+\mathbf{v}\in[\circ]$.

	We see $x+\mathbf{v}\in [\bullet]$ for all $x\in[\circ]$ and $\mathbf{v}\in V_\circ$
	in the same way.

	(3)
    Take an arbitrary $\mathbf{v}\in V_\bullet$.
    Then $\mathbf{v}=\mathbf{v}(e)$ for some $e\in\Phi(\mathcal{E})_x=\Phi(\mathcal{E}_{x'})$.
	Here, $x\in[\bullet]$ and $x'\in [\bullet]_\mathcal{X}$ satisfy $\mathbf{v}\in V_x$ and $x=\Phi(x')$.
    By definition, the inverse directed edge $e^{-1}\in\Phi(\mathcal{E})$ satisfies
    $o(e^{-1})=t(e)\in[\circ]$ and $t(e^{-1})=o(e)\in[\bullet]$.
    Hence Assertion~(1) implies $\mathbf{v}(e^{-1})\in V_\circ$.
    Moreover,
    \begin{align*}
        \mathbf{v}(e^{-1})
        =
        	t(e^{-1})-o(e^{-1})
        =
        	o(e)-t(e)
        =
	        -\mathbf{v}(e)
        =
    	    -\mathbf{v}.
    \end{align*}
    Therefore $-\mathbf{v}\in V_\circ$, which implies $V_\bullet\subset -V_\circ$.
    In a similar way, we can show $V_\bullet\supset -V_\circ$.
\end{proof}

We conclude this subsection with examples
in $\RealNum^d$ having $\mathcal{D}_m$ as a base graph. 
\begin{example}\label{ex84294820942}
	As explained in the introduction, we have the following examples.
	\begin{enumerate}
		\item\label{item890284309194}
				The $d$-dimensional square lattice in $\RealNum^d$ has $\mathcal{D}_{2d}$ as a base graph.
				We have $V_\bullet = V_\circ = \{\pm \mathbf{e}_1, \dots, \pm \mathbf{e}_d\}$.
		\item	The hexagonal lattice in $\RealNum^2$ has $\mathcal{D}_3$ as a base graph,
				and we have
				\begin{align*}
					V_\bullet
					&=
						\left\{
							+\mathbf{e}_1,
							\frac{1}{2}
							\big(
								-\mathbf{e}_1\pm \sqrt{3}\mathbf{e}_2
							\big)
						\right\},
					&
					V_\circ
					&=
						-
						V_\bullet.
				\end{align*}
		\item	The brick-wall lattice  in $\RealNum^2$ has $\mathcal{D}_3$ as a base graph,
				and we have
				$
					V_\bullet = \{\pm\mathbf{e}_1,\mathbf{e}_2\}
				$
				and
				$
					V_\circ
					=
						-
						V_\bullet
				$.
		\item	From the general theory of topological crystallography,
				we have the following.
				The first homology group of $\mathcal{D}_m$ is given by
				$H_1(\mathcal{D}_m,\Integers) \cong \Integers^{m-1}$,
				since the independent cycles are generated by differences of edges.
				Choosing a primitive subgroup
				$\Gamma \subset H_1(\mathcal{D}_m,\Integers)$ such that 
				$H_1(\mathcal{D}_m,\Integers) / \Gamma \cong \Integers^d$,
				we obtain a $\Integers^d$-periodic covering graph of $\mathcal{D}_m$,
				which can be realized as a periodic graph in $\RealNum^d$.
				For more details, see \cite{Sunada2013Book}.
	\end{enumerate}
\end{example}

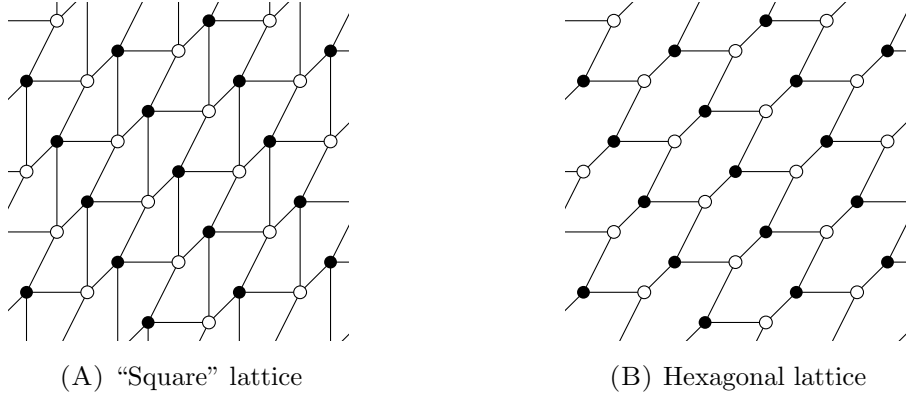
\begin{figure}
	\begin{minipage}[t]{0.45\linewidth}
		\centering
		\begin{tikzpicture}[scale=0.8]
			\clip (-2.8,-2.8) rectangle (2.8,2.8);
			
			\coordinate (vectorU) at ($0.5*(2,0)$); 
			\coordinate (vectorV) at ($0.5*(1,2)$); 
			\coordinate (vectorW) at ($0.5*(-1,-1)$);

			\foreach \i in {-5,...,5}{
				\foreach \j in {-5,...,5}{
					\coordinate (A) at ($\i*(vectorU)-\i*(vectorV)+\j*(vectorU)-\j*(vectorW)$);
					\draw (A) -- ($(A)+(vectorU)$);
					\draw (A) -- ($(A)+(vectorV)$);
					\draw (A) -- ($(A)+(vectorW)$);
					\draw (A) -- ($(A)+(vectorU)-(vectorV)+(vectorW)$);
				}
			}

			\foreach \i in {-5,...,5}{
				\foreach \j in {-5,...,5}{
					\coordinate (A) at ($\i*(vectorU)-\i*(vectorV)+\j*(vectorU)-\j*(vectorW)$);
					\fill (A) circle (3pt);
					\filldraw[fill=white] ($(A)+(vectorU)$) circle (3pt);
				}
			}
		\end{tikzpicture}
		\subcaption{``Square'' lattice}
		\label{fig433423242432342}
	\end{minipage}
	\begin{minipage}[t]{0.45\linewidth}
	    \centering
		\begin{tikzpicture}[scale=0.8]
			\clip (-2.8,-2.8) rectangle (2.8,2.8);
			
			\coordinate (vectorU) at ($0.5*(2,0)$); 
			\coordinate (vectorV) at ($0.5*(1,2)$); 
			\coordinate (vectorW) at ($0.5*(-1,-1)$);

			\foreach \i in {-5,...,5}{
				\foreach \j in {-5,...,5}{
					\coordinate (A) at ($\i*(vectorU)-\i*(vectorV)+\j*(vectorU)-\j*(vectorW)$);
					\draw (A) -- ($(A)+(vectorU)$);
					\draw (A) -- ($(A)+(vectorV)$);
					\draw (A) -- ($(A)+(vectorW)$);
				}
			}

			\foreach \i in {-5,...,5}{
				\foreach \j in {-5,...,5}{
					\coordinate (A) at ($\i*(vectorU)-\i*(vectorV)+\j*(vectorU)-\j*(vectorW)$);
					\fill (A) circle (3pt);
					\filldraw[fill=white] ($(A)+(vectorU)$) circle (3pt);
				}
			}
		\end{tikzpicture}
	    \subcaption{Hexagonal lattice}
		\label{fig490329401013423}
	\end{minipage}
	\caption{``Distorted'' lattices}
	\label{fig40239401940231}
\end{figure}

Next, we consider examples of ``distorted'' lattices.
\begin{example}
	By periodically deforming the square and hexagonal lattices,
	we obtain the following examples.
	\begin{enumerate}
		\item	For the graph in \figref{fig433423242432342}, 
				we set $V_\bullet=\{v_1,v_2,v_3,v_4\}$ and $V_\circ=-V_\bullet$ with
				\begin{align*}
					v_1 &= \mathbf{e}_1, &
					v_2 &= \frac{1}{2}(\mathbf{e}_1+2\mathbf{e}_2), &
					v_3 &= -\frac{1}{2}(\mathbf{e}_1+\mathbf{e}_2), &
					v_4 &= v_1-v_2+v_3= -\frac{3}{2}\mathbf{e}_2.
				\end{align*}
		\item	For the graph in \figref{fig490329401013423},
				we set $V_\bullet=\{v_1,v_2,v_3\}$ and $V_\circ=-V_\bullet$ with
				\begin{align*}
					v_1 &= \mathbf{e}_1, &
					v_2 &= \frac{1}{2}(\mathbf{e}_1+2\mathbf{e}_2), &
					v_3 &= -\frac{1}{2}(\mathbf{e}_1+\mathbf{e}_2).
				\end{align*}
	\end{enumerate}
\end{example}

\subsection{Elephant random walks} \label{sec9092423242092}
Let $\mathcal{S}_\even, \mathcal{S}_\odd, V_\even, V_\odd \subset \RealNum^d$.
Assume $V_\odd$ and $V_\even$ are finite.
We define the ERW $\{S_n\}_{n\geq 0}$ on $\mathcal{S}=\mathcal{S}_\even\cup \mathcal{S}_\odd$,
with increments $X_n=S_n-S_{n-1}$ ($n\geq 1$) taking values in $V = V_\odd \cup V_\even$.
More precisely, we define the ERW such that $S_n\in \mathcal{S}_\even$ and $X_{n+1}\in V_\odd$ for even $n$,
and $S_n\in \mathcal{S}_\odd$ and $X_{n+1}\in V_\even$ for odd $n$.
To this end, we introduce some assumptions.
\begin{assumption}\label{assumption49023492104}
	We assume the following.
	\begin{enumerate}
		\item\label{item100} 	$x+v\in \mathcal{S}_\odd$ for $x\in \mathcal{S}_\even$ and $v\in V_\odd$.
		\item\label{item200}	$x+v\in \mathcal{S}_\even$ for $x\in \mathcal{S}_\odd$ and $v\in V_\even$.
		\item\label{item300}	$V_\odd=-V_\even$.
		\item\label{item400} 	$\zerovector_d\in \mathcal{S}_\even$.
		\item\label{item500}	$V_\odd\cap V_\even=\emptyset$.
	\end{enumerate}
\end{assumption}
Assumptions \itemref{item100}, \itemref{item200}, and \itemref{item300}
are natural in view of \pref{prop4892384092}.
With the notation in \secref{sec487239842093841},
we have the identification
$\mathcal{S}_\even=[\bullet]$, $\mathcal{S}_\odd=[\circ]$,
$V_\odd=V_\bullet$ and $V_\even=V_\circ$.
However, $[\bullet]$, $[\circ]$, $V_\bullet$, and $V_\circ$ themselves are not used hereafter.
Assumption \itemref{item400} ensures that the ERW starts from the origin,
while Assumption \itemref{item500} restates the disjointness condition in \secref{sec940290329}.
Assumptions \itemref{item300} and \itemref{item500} imply $\zerovector_d\not\in V$. 
From Assumption \itemref{item300}, we can write $m=\sharp V_\odd = \sharp V_\even$.
Since the ERW moves in directions of $V_\odd$ and $V_\even$,
the case $m\le1$ is trivial.
In what follows, we assume that $m\ge2$.
Note that once we obtain $\{X_n\}_{n\geq 1}$, we can construct $\{S_n\}_{n\geq 0}$ as follows.
Set $S_0=\zerovector_d\in \mathcal{S}_\even$ and $S_{n+1}=S_n+X_{n+1}$.
Here, we note that, for even $n$, we have $S_n\in \mathcal{S}_\even$ and $X_{n+1}\in V_\odd$,
which implies $S_{n+1}=S_n+X_{n+1}\in\mathcal{S}_\odd$.
For odd $n$, $S_n\in \mathcal{S}_\odd$ and $X_{n+1}\in V_\even$
implies $S_{n+1}=S_n+X_{n+1}\in\mathcal{S}_\even$.
In what follows, we work under the abstract conditions stated above;
by \pref{prop4892384092}, $d$-dimensional topological crystals over dipole graphs
provide a natural class of examples satisfying these assumptions.

We now define the increments $X_n$ of the ERW in terms of $V_\odd$ and $V_\even$.
Let $X_1$ be a random variable uniformly distributed over the set $V_\odd$.
For $n\geq 1$, we recursively define $X_{n+1}$ according to the following transition rule.
Let $p,q\in[0,1]$.
The elephant first chooses a random time $\mathcal{U}_n$ uniformly from $\{1,\ldots,n\}$
and remembers the increment $X_{\mathcal{U}_n}$.
If $X_{\mathcal{U}_n}$ is admissible at time $n+1$,
the elephant repeats it with probability $p$
and otherwise chooses uniformly from the other admissible directions.
If $X_{\mathcal{U}_n}$ is not admissible, then $-X_{\mathcal{U}_n}$ is admissible by
$V_\odd=-V_\even$, and the elephant chooses $-X_{\mathcal{U}_n}$ with probability $q$
and otherwise chooses uniformly from the other admissible directions.
To formulate this rule precisely, for $v\in V$ and $k\geq2$, let
$G_k^\odd(v)$ and $G_k^\even(v)$ be random variables taking values in $V_\odd$ and $V_\even$,
respectively, whose distributions are defined as follows.
Let $w \in V_\odd$.
\begin{itemize}
	\item If $v \in V_\odd$, then we define
		\begin{align*}
			\prob(G^\odd_k(v) = w)
			=
				\begin{cases}
					p & \text{if $ w=v $}, \\
					\frac{1-p}{m-1} & \text{if $ w \neq v $}.
				\end{cases}
		\end{align*}
	\item If $v \in V_\even$ and hence $- v \in V_\odd$, then we define
		\begin{align*}
			\prob(G^\odd_k(v) = w)
			=
				\begin{cases}
					q & \text{if $w=-v$}, \\
					\frac{1-q}{m-1} & \text{if $ w\neq -v $}.
				\end{cases}
		\end{align*}
\end{itemize}
Next, we introduce the distribution of $G^\even_k(v)$ as follows.
Let $w \in V_\even$.
\begin{itemize}
	\item If $v \in V_\even$, then we define
		\begin{align*}
			\prob(G^\even_k(v) = w)
			=
				\begin{cases}
					p & \text{if $w = v$}, \\
					\frac{1-p}{m-1} & \text{if $w \neq v$}.
				\end{cases}
		\end{align*}
	\item If $v \in V_\odd$ and hence $- v \in V_\even$, then we define
		\begin{align*}
			\prob(G^\even_k(v) = w)
			=
				\begin{cases}
					q & \text{if $ w=-v $}, \\
					\frac{1-q}{m-1} & \text{if $ w \neq -v $}.
				\end{cases}
		\end{align*}
\end{itemize}
Set 
\begin{align*}
	G_k(v)
	=
		G^\odd_k(v) \indicator{\odd}(k)
		+
		G^\even_k(v) \indicator{\even}(k).
\end{align*}
Then $G_k(v)$ is $V_\odd$-valued if $k$ is odd and $V_\even$-valued if $k$ is even.
With this preparation, we define, for $n\geq 1$,
\begin{align*}
	X_{n+1} = G_{n+1}(X_{\mathcal{U}_n}).
\end{align*}
Assume that $X_1, \mathcal{U}_1, \{ G_2 (v) \}_{v\in V}, \mathcal{U}_2, \{ G_3 (v) \}_{v\in V},\dots$ are independent.
Denote by $\sigmaField_0$ the trivial $\sigma$-field
and write $\sigmaField_1= \sigma(X_1)$ and, for $n\geq 2$, 
\begin{align*}
	\sigmaField_n = \sigma(X_1, \mathcal{U}_1, \{G_2(v)\}_{v \in V},\dots, X_{n-1}, \mathcal{U}_{n-1}, \{G_n(v)\}_{v \in V}).
\end{align*}
Note that $X_n$ is $\sigmaField_n$-measurable.
We define the position process of the ERW $\{S_n\}_{n\ge0}$ by $S_0 = 0$ and, for $n\geq 0$,
\begin{align*}
	S_{n+1} = S_n+X_{n+1}.
\end{align*}

Here, we make a remark on $G^\odd_k(v)$ and $G^\even_k(v)$
and introduce important quantities to analyze $\{S_n\}_{n\ge0}$.
From the definition above, we have
\begin{equation}\label{eq:Gk-expectation}
\begin{aligned}
  \expect[G^\odd_k(v)]
  &=
    \begin{cases}
      \alpha v + (1 - \alpha) \bar{v} & \text{if $v \in V_\odd$}, \\
      -\beta  v + (1 - \beta) \bar{v} & \text{if $v \in V_\even$},
    \end{cases}
  \\
  \expect[G^\even_k(v)]
  &=
    \begin{cases}
      \alpha v - (1 - \alpha) \bar{v} & \text{if $v \in V_\even$}, \\
      -\beta  v - (1 - \beta) \bar{v} & \text{if $v \in V_\odd$}.
    \end{cases}
\end{aligned}
\end{equation}
Here 
\begin{align*}
    \alpha
    &=
        \frac{mp-1}{m-1},
    &
    \beta
    &=
        \frac{mq-1}{m-1},
    &
    \bar{v}
    &=
        \frac{1}{m}\sum_{v \in V_\odd} v.
\end{align*}
Since $V_\odd = - V_\even$, we have $\bar{v} = - \frac{1}{m}\sum_{v \in V_\even} v$.
For later use, it is convenient to rewrite the expressions in terms of the sum and difference of $\alpha$ and $\beta$,
that is, we set
\begin{align*}
    \gamma
    &=
        \frac{1}{2}(\alpha + \beta),
    &
    \delta
    &=
        \frac{1}{2}(\alpha - \beta).
\end{align*}
Since $-\frac{1}{m-1}\leq \alpha,\beta\leq 1$,
the parameters $\gamma$ and $\delta$ run over
\begin{align}
	\label{eq4539028429084}
	-
	\frac{1}{m-1}
	&\leq
		\gamma+\delta
	\leq
		1,
	&
	-
	\frac{1}{m-1}
	&\leq
		\gamma-\delta
	\leq
		1.
\end{align}
In particular, we have $-1\leq \delta\leq 1$ for $m=2$.
We obtain the following results for the ERW
in the case $\gamma=1$ or $\delta=1$.
\begin{remark}\label{rem:degenerate-case}
	The following discussion shows that
	the ERW becomes deterministic after the initial step
	if $\gamma=1$ or $\delta=1$.
	Here we recall $V_\odd\cap V_\even=\emptyset$.
	\begin{itemize}
		\item	If $\gamma = 1$, then we have $\alpha = \beta = 1$,
				that is, $p=q=1$.
				In this case, we have $X_2=-X_1$ from the definition.
				Thereafter, we see that the increment of the elephant
				alternates between $X_1$ and $-X_1$ forever,
				that is, $X_n = (-1)^{n-1} X_1$ for all $n$.
		\item	Assume that $\delta = 1$.
				Note that $\delta=1$ holds if and only if $m=2$, $p=1$ and $q=0$.
				Let $V_\odd = \{v_1, v_2\}$ and $V_\even = \{-v_1, -v_2\}$.
				Here, we can assume that $v_1\neq \zerovector_d$, $v_2\neq \zerovector_d$, $v_1 \neq v_2$ and $v_1 \neq -v_2$.
				Then, the dynamics becomes deterministic once the initial step is fixed.
				Indeed, for $v, w\in V_\odd$ with $v\neq w$, we have
				\begin{align*}
					\prob( X_{2k-1}=v,\ X_{2k}=-w , \forall k\ge1 )=\prob( X_1=v ).
				\end{align*}
				Hence the walk alternates between two fixed directions, and the only randomness
				comes from $X_1$.
	\end{itemize}
	Therefore, we assume that $\gamma<1$ and $\delta<1$ in what follows.
\end{remark}

We finally introduce several processes associated with $\{S_n\}_{n\geq 0}$.
For each $v\in V$, let $Y_n(v)$ denote the number of times that $v$ has been chosen up to time $n$.
Then, it follows that $Y_0(v)=0$ and, for $n\geq 0$, 
\begin{align}\label{eq:Yn-recurrence}
	Y_{n+1}(v)
	=
		Y_n(v) + Z_{n+1}(v),
\end{align}
where
\begin{align*}
	Z_{n+1}(v)=\setindicator{X_{n+1} = v}.
\end{align*}
For notational simplicity, we write $Y_n=(Y_n(v))_{v\in V}$ and $Z_n=(Z_n(v))_{v\in V}$
and treat them as vectors in $\RealNum^{2m}$. Then, we can write \eqref{eq:Yn-recurrence} as $Y_{n+1}=Y_n + Z_{n+1}$.
From the definition, we have
\begin{align}\label{eq48902890218401984}
	X_k
	=
		\sum_{v\in V} v\setindicator{X_k = v}
	=
		\sum_{v\in V} vZ_k(v)
\end{align}
From this and $Y_n(v)=\sum_{k=1}^n Z_k(v)$, we have
\begin{align}
	\label{eq453849023809481}
	S_n
	=
		\sum_{k=1}^n
			X_k
	=
		\sum_{k=1}^n
		\sum_{v\in V}
			vZ_k(v)
	=
		\sum_{v \in V} v Y_n(v).
\end{align}
Furthermore, we introduce the auxiliary process
\begin{align}
	\label{eq843904820941423}
	T_n 
	&=
		\sum_{v \in V_\odd} v Y_n(v) 
		-
		\sum_{v \in V_\even} v Y_n(v).
\end{align}

\subsection{Main results}\label{sec5409480923184092}
We are in a position to state our main results.
From \rref{rem:degenerate-case}, it is natural to assume $\gamma<1$ and $\delta < 1$.
As stated in the introduction, we assume $V_\odd\cap V_\even=\emptyset$.
The first two results concern $Y_n$, and 
their proofs are based on the P{\'o}lya urn approach.
\begin{theorem}\label{thm43092043241902494}
	For all $n\in\NaturalNum$ and $v\in V$, we have
	\begin{align*}
		\expect[Y_n(v)]
		=
			\begin{cases}
				\frac{1}{m}
				\lceil{\frac{n}{2}}\rceil & \text{if $v\in V_\odd$,}\\
				\frac{1}{m}
				\lfloor{\frac{n}{2}}\rfloor  & \text{if $v\in V_\even$.}
			\end{cases}
	\end{align*}
	Here, $\lceil\cdot\rceil$ and $\lfloor\cdot\rfloor$ denote the ceiling and floor functions, respectively.
	In particular, we have
	\begin{align*}
		\expect[S_n]
		=
		    \begin{cases}
				\bar{v} & \text{if $n$ is odd}, \\
				\zerovector_d & \text{if $n$ is even}.
			\end{cases}
	\end{align*}
	Hence, for all $v\in V$, we have
	\begin{align*}
		\lim_{n\to\infty} \frac{\expect[Y_n(v)]}{n} 
		&=
			\frac{1}{2m},
		&
		\lim_{n\to\infty} \frac{\expect[S_n]}{n} 
		&=
			\zerovector_d.
	\end{align*}
\end{theorem}
\begin{theorem}\label{thm:SLLN}
	Assume that $\gamma<1$ and $\delta < 1$.
	Then, for all $v\in V$, we have
	\begin{align*}
		\lim_{n\to\infty} \frac{Y_n(v)}{n} = \frac{1}{2m} \quad \text{a.s.}
	\end{align*}
	In particular, we have
	\begin{align*}
		\lim_{n\to\infty} \frac{S_n}{n} = \zerovector_d \quad \text{a.s.}
	\end{align*}
\end{theorem}
The first assertion in \tref{thm:SLLN} states that the elephant chooses each vector in $V$
asymptotically with the same frequency $1/(2m)$ almost surely.
Using this result and \eqref{eq453849023809481}, we immediately obtain the strong law
of large numbers for the position process $S_n$.
Proofs of \tref[thm43092043241902494]{thm:SLLN} are given in \secref{sec320492094209410}.

Previous studies on the ERW showed that the model exhibits phase transitions
depending on the memory parameters.
To characterize the phase transition, 
we describe the asymptotic behavior of the centered second moment of $S_n$.
\begin{theorem}\label{thm:Asymptotic-second-moment}
	Assume that $\gamma<1$ and $\delta < 1$.
	Then, the following limit exists and is positive:
	\begin{align}\label{eq489208401}
		\sigma^2 
		=
			\lim_{n\to\infty} \expect[|X_n - \expect[X_n]|^2].
	\end{align}
	In addition, there exists a constant $C_{\alpha,\beta}>0$
	depending on $\alpha$ and $\beta$ such that, as $n\to\infty$,
	\begin{align*} 
		\expect[|S_n - \expect[S_n]|^2]
		\sim
			\begin{cases}
				\frac{\sigma^2}{1-2\delta}n & \text{if $\delta<\frac{1}{2}$}, \\
				\sigma^2 \, n\log n & \text{if $\delta=\frac{1}{2}$}, \\
				C_{\alpha, \beta} n^{2\delta}  & \text{if $\frac{1}{2}<\delta<1$}.
			\end{cases}
	\end{align*}
\end{theorem}
The first assertion of this theorem is a corollary of \tref{thm:SLLN}
and is proved in \secref{sec320492094209410}.
The second assertion, proved in \secref{sec:Centered-second-moment},
shows that the model exhibits three different regimes, which we define below.
\begin{definition}
	We say that the ERW is diffusive if $\delta < 1/2$,
	critical if $\delta = 1/2$, and superdiffusive if $\delta > 1/2$.
\end{definition}
This classification of regimes is consistent with previous studies on the ERW.
For example, the MERW in Bercu and Laulin 
\cite{BercuLaulin2019MERW} exhibits the same phase transition
depending on the memory parameter $\alpha = (2dp-1)/(2d-1)$.
Therefore, we analyze the precise limit theorems in each of the three regimes.

Next we state our main results in each of the three regimes.
We first record the following consequences of
\eqref{eq4539028429084} and \rref{rem:degenerate-case}.
\begin{itemize}
	\item	In the diffusive regime $\delta < 1/2$,
			although $\gamma=1$ may occur,
			we exclude this case.
	\item	In the critical regime $\delta = 1/2$, we have $\gamma<1$.
			Hence any assumption on $\gamma$ is unnecessary.
	\item	In the superdiffusive regime $\delta > 1/2$,
			any assumption on $\gamma$ is unnecessary
			while the case $\delta=1$ is excluded.
\end{itemize}
We are now in a position to state limit theorems in each regime.

In the diffusive regime $\delta < 1/2$, 
under $\gamma<1$, we obtain 
the strong law of large numbers with a non-standard normalization 
and the central limit theorem.
\begin{theorem}\label{thm:SLLN-in-diffusive}
	Assume that $\gamma<1$ and $\delta < 1/2$.
	Then, for any $\eta>0$, we have
	\begin{align*}
		|S_n|
		= 
			o(\sqrt{n(\log n)^{1+\eta}})
		\quad
		\text{a.s. as\,\, $n\to\infty$.}
	\end{align*}
\end{theorem}
\begin{theorem}\label{thm:CLT-in-diffusive}
	Assume that $\gamma<1$ and $\delta < 1/2$.
	Then, we have
	\begin{align*}
		\lim_{n\to\infty}
			\frac{S_n}{\sqrt{n}}
		=
			\mathcal{N}_d
				\left(
					\zerovector_d,
					\frac{1}{m(1-2\delta)}
					\sum_{v\in V_\odd}
						(v - \bar{v})(v - \bar{v})^\top
				\right)
		\quad \text{in distribution.}
	\end{align*}
\end{theorem}

In the critical regime $\delta = 1/2$,
the following results hold without any assumption on $\gamma$.
\begin{theorem}\label{thm:SLLN-in-critical}
	Assume that $\delta=1/2$.
	Then, for any $\eta>0$, we have
	\begin{align*}
		|S_n|
		= 
			o(\sqrt{n (\log n)(\log \log n)^{1+\eta}})
		\quad
		\text{a.s. as\,\, $n\to\infty$.}
	\end{align*}
\end{theorem}
\begin{theorem}\label{thm:CLT-in-critical}
	Assume that $\delta=1/2$.
	Then, we have
	\begin{align*}
		\lim_{n\to\infty}
			\frac{S_n}{\sqrt{n \log n}}
		=
			\mathcal{N}_d
				\left(
					\zerovector_d,
				 	\frac{1}{m}
					\sum_{v\in V_\odd}
						(v - \bar{v})(v - \bar{v})^\top
				\right)
		\quad \text{in distribution.}
	\end{align*}
\end{theorem}

In the superdiffusive regime $\delta > 1/2$,
we obtain the strong law of large numbers with a non-standard normalization.
\begin{theorem}\label{thm:SLLN-in-superdiffusive}
	Assume that $1/2<\delta<1$.
	Then, there exists a random vector $L \in \RealNum^d$ such that
	\begin{align*}
		\lim_{n\to\infty} \frac{S_n}{n^\delta} = L \quad \text{a.s. and\,  in $L^2$.}
	\end{align*}
	Furthermore, we have
	\begin{align*}
		\expect[L]
		&=
			\zerovector_d,
		&
		\expect[|L|^2]
		&=
			C_{\alpha,\beta},
	\end{align*}
	where $C_{\alpha,\beta}$ is a positive constant appearing in \tref{thm:Asymptotic-second-moment}.
	In particular, $L$ is not constant.
\end{theorem}
Note that constants $\sigma^2$, $C_{\alpha,\beta}$ in \tref{thm:Asymptotic-second-moment},
$
	\frac{1}{m}
	\sum_{v\in V_\odd}
		(v - \bar{v})
		(v - \bar{v})^\top
$
in \tref[thm:CLT-in-diffusive]{thm:CLT-in-critical}
depend on $V_\odd$ and $V_\even$.
In graph settings, this means that 
they depend on the geometric structure of the graphs.
Note that the matrix 
$
	\frac{1}{m}
	\sum_{v\in V_\odd}
		(v - \bar{v})
		(v - \bar{v})^\top
$
may be singular under the abstract assumptions.
For example, $V_\odd=\{\mathbf{e}_1,\mathbf{e}_2\}$ 
and $V_\even=\{-\mathbf{e}_1,-\mathbf{e}_2\}$.
The above theorems are shown in \secref{sec:limit-theorems}.

\section{The urn approach}\label{sec320492094209410}
In this section, we give an expression for 
$\expect[Z_{n+1}|\sigmaField_n]$ and show \tref[thm43092043241902494]{thm:SLLN}.
Recall that $V_\odd\cap V_\even=\emptyset$.
We now index the elements of $V_\odd$, $V_\even$, and $V=V_\odd\cup V_\even$
and introduce matrix notation for $Y_n$ and $Z_n$.
Taking into account the relation $V_\odd=-V_\even$,
we index the elements such that $v_i=-v_{m+i}$ for all $1\leq i\leq m$ and write
\begin{align*}
	V_\odd
	&=
		\{v_1,\dots,v_m\},
	&
	V_\even
	&=
		\{v_{m+1},\dots,v_{2m}\}.
\end{align*}
We also introduce a matrix $\mathbf{V}\in \RealNum^{d\times 2m}$ by
\begin{align*}
	\mathbf{V}
	=
		(
			v_1, \dots,v_m,
			v_{m+1},\dots,v_{2m}
		).	
\end{align*}
Using the indexing above,
we also introduce indices for $Y_n$ and $Z_n$ by writing
$
	Y_n 
	=
		(Y_n^1,\dots,Y_n^{2m})^\top
$
with $Y_n^j=Y_n(v_j)$
and
$
	Z_n 
	=
		(Z_n^1,\dots,Z_n^{2m})^\top
$
with  $Z_n^j=Z_n(v_j)$.
From \eqref{eq48902890218401984} and \eqref{eq453849023809481},
it is easy to see that
\begin{align*}
	X_{n+1}
	&=
		\sum_{i=1}^{2m}
			v_i Z_{n+1}^i
	=
		\mathbf{V} Z_{n+1}	,
	&
	S_n
	&=
		\sum_{i=1}^{2m}
			v_i Y_n^i
	=
		\mathbf{V} Y_n.
\end{align*}

Hereafter, we confirm that $\{Y_n\}_{n\geq 0}$ forms an urn process
that fits into the framework of Yan-Cheng-Bai \cite{YanChengBai2006}.
We first rewrite the recurrence relation \eqref{eq:Yn-recurrence} of $Y_n$.
To this end, we introduce a matrix-valued random variable
$D_n=(D_n^{i,j})_{1\leq i,j\leq 2m}$ by 
$
	D_n^{i,j}
	=
		\setindicator{G_n(v_j) = v_i}
$.
As $D_n$ is independent of $\sigmaField_{n-1}$, 
we have $\expect[D_n|\sigmaField_{n-1}]=\expect[D_n]$.
Hence, we set $H_n =\expect[D_n]$ for $n\geq 2$.
From the definition of $Z_n^i$, we have
\begin{align*}
    Z_{n+1}^i
    =
        \setindicator{X_{n+1} = v_i}
    =
        \sum_{j=1}^{2m}
            \setindicator{G_{n+1}(v_j) = v_i} \setindicator{X_{\mathcal{U}_n} = v_j}
    =
        \sum_{j=1}^{2m}
            D_{n+1}^{i,j} Z_{\mathcal{U}_n}^j
\end{align*}
and we can rewrite \eqref{eq:Yn-recurrence} as $Y_{n+1} = Y_n + D_{n+1} Z_{\mathcal{U}_n}$.
We also have
\begin{align}\label{eq:Conditional-expectation-of-urnincrement}
    \expect[Z_{n+1}|\sigmaField_n]
    =
        \sum_{k=1}^{n} \prob(\mathcal{U}_n = k) \expect[D_{n+1} Z_k |\sigmaField_n] 
    =
        \frac{\expect[D_{n+1}]}{n} \sum_{k=1}^{n} Z_k 
    =
        \frac{H_{n+1}}{n} Y_n.
\end{align}
Thus, $\{Y_n\}_{n\ge1}$ is a $2m$-color P\'olya-type urn process with replacement matrix $D_{n+1}$.

Next, for $n\geq 2$, we see
\begin{align*}
	H_n
	=
		\indicator{\odd}(n)
		H_\odd
		+
		\indicator{\even}(n)
		H_\even.
\end{align*}
Here, using 
\begin{align*}
	A
	&=
		\alpha I_m+\frac{1-\alpha}{m}\onevector_m\onevector_m^\top,
	&
	B
	&=
		\beta I_m+\frac{1-\beta}{m}\onevector_m\onevector_m^\top,
\end{align*}
and the $2\times 2$ block matrix notation, we set
\begin{align*}
	H_\even
	&=
		\begin{pmatrix}
			O & O\\
			B & A\\
		\end{pmatrix},
	&
	H_\odd
	&=
		\begin{pmatrix}
			A & B\\
			O & O
		\end{pmatrix}.
\end{align*}
Indeed, noting that 
\begin{align*}
	D_n^{i,j}
	=
		\setindicator{G_n^\odd(v_j)\indicator{\odd}(n)+G_n^\even(v_j)\indicator{\even}(n)=v_i}
	=
		\setindicator{G_n^\odd(v_j)=v_i}
		\indicator{\odd}(n)
		+
		\setindicator{G_n^\even(v_j)=v_i}
		\indicator{\even}(n),
\end{align*}
and recalling that $G_n^\odd$ and $G_2^\odd$ have the same distribution,
as do $G_n^\even$ and $G_2^\even$,
we have
\begin{align*}
	H_n^{i,j}
	=
		\expect[D_n^{i,j}]
	=
		\prob(G_2^\odd(v_j)=v_i)
		\indicator{\odd}(n)
		+
		\prob(G_2^\even(v_j)=v_i)
		\indicator{\even}(n).
\end{align*}
Hence, we obtain the stated expressions for $H_n$.

We are in a position to show \tref{thm43092043241902494}.

\begin{proof}[Proof of \tref{thm43092043241902494}]
	The main step of the proof is to show that
	there exist sequences $(a_n)$ and $(b_n)$ such that
	\begin{align*}
		\expect[Y_n]
		=
			\begin{pmatrix}
				a_n \onevector_m\\
				b_n \onevector_m\\
			\end{pmatrix}
	\end{align*}
	and to derive their recurrence relations.
	We proceed by induction.
	From the definition, by choosing $a_0=0$, $b_0=0$, $a_1=\frac{1}{m}$ and $b_1=0$,
	we obtain the assertion for $n=0,1$.
	Assume that the assertion holds for some $n\geq 1$.
	Then, by \eqref{eq:Conditional-expectation-of-urnincrement}, we have
	\begin{align*}
		\expect[Y_{n+1}]
		=
			\left\{
				I_{2m}
				+
				\frac{H_{n+1}}{n}
			\right\}
			\expect[Y_n]
		=
			\left\{
				I_{2m}
				+
				\frac{1}{n}
				\left(
					\indicator{\odd}(n+1)H_\odd
					+
					\indicator{\even}(n+1)H_\even
				\right)
			\right\}
			\begin{pmatrix}
				a_n \onevector_m\\
				b_n \onevector_m\\
			\end{pmatrix}.
	\end{align*}
	Here we have
	\begin{align*}
		H_\odd
		\begin{pmatrix}
			a_n \onevector_m\\
			b_n \onevector_m\\
		\end{pmatrix}
		=
			\begin{pmatrix}
				A & B\\
				O & O
			\end{pmatrix}
			\begin{pmatrix}
				a_n \onevector_m\\
				b_n \onevector_m\\
			\end{pmatrix}
		=
			\begin{pmatrix}
				a_n A \onevector_m + b_n B \onevector_m\\
				\zerovector_m
			\end{pmatrix}
		=
			\begin{pmatrix}
				(a_n + b_n)\onevector_m\\
				\zerovector_m
			\end{pmatrix}.
	\end{align*}
	A similar result for $H_\even$ also holds.
	In addition, we see that $\onevector_{2m}^\top Y_n=n$ yields 
	$n=\expect[\onevector_{2m}^\top Y_n]=m(a_n+b_n)$.
	Hence we obtain
	\begin{align*}
		\expect[Y_{n+1}]
		=
			\begin{pmatrix}
				\left\{
					a_n + \frac{1}{m} \indicator{\odd} (n+1)
				\right\}
				\onevector_m\\
				\left\{
					b_n + \frac{1}{m} \indicator{\even} (n+1)
				\right\}
				\onevector_m
			\end{pmatrix},
	\end{align*}
	which implies
	\begin{align*}
		a_{n+1} &= a_n + \frac{1}{m} \indicator{\odd} (n+1), 
		&
		b_{n+1} &= b_n + \frac{1}{m} \indicator{\even} (n+1).
	\end{align*}
	The main step is thus established.

	From the main step, we have
	$a_n=\frac{1}{m}\lceil{\frac{n}{2}}\rceil$ and $b_n=\frac{1}{m}\lfloor{\frac{n}{2}}\rfloor$.
	Since $S_n = \mathbf{V} Y_n$, we have
	\begin{align*}
		\expect[S_n]
		=
		    \mathbf{V} \expect[Y_n]
		=
		    a_n \sum_{i=1}^{m} v_i + b_n \sum_{i=m+1}^{2m} v_i 
		=
			(a_n - b_n) m\bar{v}.
	\end{align*}
	This completes the proof.
\end{proof}

To show \tref{thm:SLLN}, we introduce
\begin{align*}
	H
	=
		\frac{1}{2}\{H_\odd+H_\even\}
	=
		\frac{1}{2}
		\begin{pmatrix}
			A & B\\
			B & A
		\end{pmatrix}.
\end{align*}
Before proving \tref{thm:SLLN}, we give details on the relationship 
between \tref{thm:SLLN} and \cite[Theorem~1]{YanChengBai2006}.
\begin{remark}\label{rem489208409238091483921}
	In \cite[Theorem~1]{YanChengBai2006}, the authors used 
	the condition that all eigenvalues other than $1$ have modulus strictly less than $1$.
	However, $H$ does not satisfy this condition in general.
	Consider the case $\gamma=-1$.
	Then, $\alpha=-1$ and $\beta=-1$, so $m=2$.
	In this case, for $k\geq 1$, we have
	$
		H^{2k-1}
		=
			\frac{1}{2}
			\begin{pmatrix}
				J & J\\
				J & J
			\end{pmatrix}
	$
	and
	$
		H^{2k}
		=
			\frac{1}{2}
			\begin{pmatrix}
				I_2 & I_2\\
				I_2 & I_2
			\end{pmatrix}
	$,
	where
	$
		J
		=
			\begin{pmatrix}
				0 & 1 \\
				1 & 0 \\
			\end{pmatrix}
	$.
	Hence, $H$ is nonnegative and irreducible.
	We also observe that the eigenvalues of $H$ are $0$, $1$ and $-1$.
	Therefore, \cite[Theorem~1]{YanChengBai2006} cannot be applied directly to our setting.
	The case $\delta=-1$ is analogous.
	As we will see in \lref{lem:Diagonalization-of-H},
	the eigenvalues of $H$ are $1$, $\gamma$, $0$, and $\delta$.
	Hence, if $|\gamma|<1$ and $|\delta|<1$
	then we can use \cite[Theorem~1]{YanChengBai2006} to show \tref{thm:SLLN}.
\end{remark}
For the above reasons, and for the sake of a self-contained presentation and the reader's convenience,
we give a direct proof covering the case $\gamma=-1$ or $\delta=-1$.
We divide the proof into a series of lemmas.
We first show that $H$ is diagonalizable.
\begin{lemma}\label{lem:Diagonalization-of-H}
Let $u = \frac{1}{\sqrt{m}}\onevector_m$.
For $r=1,2,\dots,m-1$, define
\begin{align*}
    w_r
    =
    	\frac{1}{\sqrt{r(r+1)}}
		\begin{pmatrix}
			\onevector_r\\
			-r\\
			\zerovector_{m-r-1}
		\end{pmatrix}
    \in\RealNum^m.
\end{align*} 
Set
\begin{align*}
	Q
	&=
		(u,w_1,\dots,w_{m-1})\in\RealNum^{m\times m},
	&
	P 
	&=
		\frac{1}{\sqrt2}
		\begin{pmatrix}
			Q & Q\\
			Q & -Q
		\end{pmatrix}
		\in\RealNum^{2m\times 2m},
\end{align*}
Then, the following holds.
\begin{enumerate}
	\item	$Q$ is an orthogonal matrix and diagonalizes $A$ and $B$.
			Explicitly, we have
			\begin{align*}
				Q^\top A Q
				&=
					\diag(1,\alpha \onevector_{m-1}^\top),
				&
				Q^\top B Q
				=
					\diag(1,\beta \onevector_{m-1}^\top).
			\end{align*}
	\item	$P$ is an orthogonal matrix and we have
			\begin{align*}
				P^\top H_\odd P
				&=
					\begin{pmatrix}
						\diag(1,\gamma \onevector_{m-1}^\top) & \diag(0,\delta \onevector_{m-1}^\top)\\
						\diag(1,\gamma \onevector_{m-1}^\top) & \diag(0,\delta \onevector_{m-1}^\top)
					\end{pmatrix},\\
				P^\top H_\even P
				&=
					\begin{pmatrix}
						\diag(1,\gamma \onevector_{m-1}^\top) & -\diag(0,\delta \onevector_{m-1}^\top)\\
						-\diag(1,\gamma \onevector_{m-1}^\top) & \diag(0,\delta \onevector_{m-1}^\top)
					\end{pmatrix},\\
				P^\top H P
				&=
					\diag(1,\gamma \onevector_{m-1}^\top,\,0,\delta \onevector_{m-1}^\top).
			\end{align*}
			Namely, $P$ diagonalizes $H$.
\end{enumerate}

\end{lemma}

\begin{proof}
	(1) Assume that $\alpha<1$.
	Since $A$ has an eigenvector $\onevector_m$ for the eigenvalue $1$ (multiplicity $1$),
	and every vector $w$ with $\onevector_m^\top w=0$ is an eigenvector for the eigenvalue $\alpha$ (multiplicity $m-1$),
	we can construct the orthogonal matrix $Q$ as above.
	In the case $\alpha=1$, $A=I_m$ has only the eigenvalue $1$ with multiplicity $m$,
	and we can also construct the orthogonal matrix $Q$ in the same way.
	Similarly, $Q$ diagonalizes $B$.

	(2) Since $Q$ diagonalizes both $A$ and $B$, the assertion holds.
\end{proof}

Next we give estimates of the following matrices:
\begin{align*}
	\Pi_{k,l}
	=
		\begin{cases}
			\left(I_{2m} + \frac{H_{k+1}}{k}\right) \cdots \left(I_{2m} + \frac{H_{l+1}}{l}\right)
			&
			\text{if $k\geq l$,}\\
			I_{2m}
			&
			\text{if $k<l$.}
		\end{cases}
\end{align*}
To this end, given a matrix $M$, we denote by $M_{i_1:i_2,j_1:j_2}$
the submatrix of $M$ consisting of the rows from $i_1$ to $i_2$ and the columns from $j_1$ to $j_2$.
If $i_1 = i_2$ (or $j_1 = j_2$), we write $M_{i_1,j_1:j_2}$ (or $M_{i_1:i_2,j_1}$) for simplicity.
For $c\geq 0$ and $k\geq l\geq 1$, we set
\begin{align*}
	\lambda_{k,l}(c)
	=
		\prod_{i=l}^{k}
			\frac{2i-1 + 2c}{2i-1}
	=
		\frac{\Gamma\left(k + \frac{1}{2} + c\right) }{\Gamma\left(k + \frac{1}{2}\right) }
		\frac{\Gamma\left(l - \frac{1}{2}\right)}{\Gamma\left(l - \frac{1}{2} + c\right)}.
\end{align*}
Write $\rho=\max\{0,\gamma,\delta\}$.
\begin{lemma}\label{lem89453204920492013}
	For $k\geq l \geq 1$, we have
	\begin{align*}
		\big\| \big(P^\top \Pi_{2k,2l-1} P\big)_{2:2m,2:2m} \big\|_\mathrm{op}
		&\le
			C
			\lambda_{k,l}(\rho),
		&
		\big\| \big(P^\top \Pi_{2k,2l} P\big)_{2:2m,2:2m} \big\|_\mathrm{op}
		&\le
			C
			\lambda_{k,l}(\rho).
	\end{align*}
	Here $C$ is a positive constant.
\end{lemma}

\begin{proof}
	The case $k=1$ is trivial.
	We show the first estimate for $k\geq 2$.
	Note
	\begin{align*}
		\Pi_{2k,2k-1}
		=
			\left(I_{2m} + \frac{H_{2k+1}}{2k}\right)
			\left(I_{2m} + \frac{H_{2k}}{2k-1}\right)
		=
			I_{2m} + \frac{2H}{2k-1} + \frac{H_\odd (H_\even -I_{2m})}{2k(2k-1)}.
	\end{align*}
	Then, \lref{lem:Diagonalization-of-H} yields
	\begin{align*}
		P^\top \Pi_{2k,2k-1} P
		&=
			\begin{pmatrix}
				1+\frac{2}{2k-1} & \zerovector_{2m-1}^\top \\
				\zerovector_{2m-1} & F_k+G_k
			\end{pmatrix},
	\end{align*}
	where
	\begin{align*}
		F_k
		&=
			I_{2m-1}
			+
			\frac{2}{2k-1}
			\diag(\gamma \onevector_{m-1}^\top,0,\delta \onevector_{m-1}^\top),\\
		G_k
		&=
			\frac{1}{2k(2k-1)}
			\begin{pmatrix}
				\gamma(\gamma-\delta-1)I_{m-1} 
					& \zerovector_{m-1}
					& \delta(\delta-\gamma-1) I_{m-1}\\
				\zerovector_{m-1}^\top
					& 0
					& \zerovector_{m-1}^\top \\
				\gamma(\gamma-\delta-1)I_{m-1} 
					& \zerovector_{m-1}
					& \delta(\delta-\gamma-1) I_{m-1} \\
			\end{pmatrix}.
	\end{align*}
	Since 
	$
		P^\top \Pi_{2k,2l-1} P
		=
			\left(
				P^\top \Pi_{2k,2k-1} P
			\right)
			\cdots
			\left(
				P^\top \Pi_{2l,2l-1} P
			\right)
	$,
	we have
	\begin{align*}
		\left(P^\top \Pi_{2k,2l-1} P\right)_{2:2m,2:2m}
		=
			(F_k+G_k)
			\cdots
			(F_l+G_l).
	\end{align*}

	We consider $\left(P^\top \Pi_{2k,2l-1} P\right)_{2:2m,2:2m}$ for $l\geq 2$.
	The above implies
	\begin{align*}
		\big\| \big(P^\top \Pi_{2k,2l-1} P\big)_{2:2m,2:2m} \big\|_\mathrm{op}
		&\leq
			\prod_{i=l}^k 
				\left\{
					\|F_i\|_\op+\|G_i\|_\op
				\right\}\\
		&=
			\sum_{S\subset\{l,\dots,k\}}
				\Bigg\{
					\prod_{i\in \{l,\dots,k\}\setminus S}
					\|F_i\|_\op
				\Bigg\}
				\Bigg\{
					\prod_{i\in S}
						\|G_i\|_\op
				\Bigg\}.
	\end{align*}
	Eigenvalues of $F_i$ are $1$, $1+\frac{2\gamma}{2i-1}$, and $1+\frac{2\delta}{2i-1}$
	and all of them are positive for $i\geq 2$.
	Hence, $\|F_i\|_\op\leq 1+\frac{2\rho}{2i-1}$, which implies
	\begin{align*}
		\prod_{i\in \{l,\dots,k\}\setminus S}
			\|F_i\|_\op
		\leq
			\prod_{i\in \{l,\dots,k\}}
				\left(
					1+\frac{2\rho}{2i-1}
				\right)
		=
			\lambda_{k,l}(\rho).
	\end{align*}
	It is clear that $\|G_i\|_\op\leq \frac{C}{i^2}$, which implies
	\begin{align*}
		\sum_{S\subset\{l,\dots,k\}}
			\prod_{i\in S}
				\|G_i\|_\op
		\leq
			\sum_{S\subset\{l,\dots,k\}}
				\prod_{i\in S}
					\frac{C}{i^2}
		=
			\prod_{i\in \{l,\dots,k\}}
				\left(
					1
					+
					\frac{C}{i^2}
				\right)
		\leq
			\exp
				\left(
					\sum_{i=1}^\infty
						\frac{C}{i^2}
				\right)
		<
			\infty.
	\end{align*}
	Hence we obtain the estimate for $P^\top \Pi_{2k,2l-1} P$ in the case $l\geq 2$.
	Since
	\begin{align*}
		\big\| \big(P^\top \Pi_{2k,1} P\big)_{2:2m,2:2m} \big\|_\mathrm{op}
		\leq
			\big\| \big(P^\top \Pi_{2k,2\cdot 2-1} P\big)_{2:2m,2:2m} \big\|_\mathrm{op}
			\big\| \big(P^\top \Pi_{2,1} P\big)_{2:2m,2:2m} \big\|_\mathrm{op}
	\end{align*}
	and $\lambda_{k,2}$ and $\lambda_{k,1}$ are comparable,
	we see the assertion for $l=1$.
	
	The second assertion for $k>l$ follows from the comparability of $\lambda_{k,l+1}$ and $\lambda_{k,l}$
	and
	\begin{align*}
		\big\| \big(P^\top \Pi_{2k,2l} P\big)_{2:2m,2:2m} \big\|_\mathrm{op}
		\leq
			\big\| \big(P^\top \Pi_{2k,2(l+1)-1} P\big)_{2:2m,2:2m} \big\|_\mathrm{op}
			\big\| \big(P^\top \Pi_{2l,2l} P\big)_{2:2m,2:2m} \big\|_\mathrm{op}.
	\end{align*}
	The case $k=l$ is trivial.
	The proof is complete.
\end{proof}	

\begin{lemma}\label{lem:Pi-trace-inequality}
	Let $E$ be a $(2m)$-dimensional nonnegative-definite and symmetric matrix with the first row and first column being zero.
	Then, for $k\geq l\geq 1$ and $i=2l-1,2l$, we have
	\begin{align*}
		\trace \left((P^\top\Pi_{2k,i} P) E (P^\top\Pi_{2k,i} P)^\top \right)
		\le
			C
			\lambda_{k,l}(\rho)^2
            \trace \left(E\right).
	\end{align*}
	Here, $C$ is a positive constant.
\end{lemma}
\begin{proof}
	For notational simplicity, we write $A=P^\top\Pi_{2k,i} P$.
	Then, since entries of the first row of $A$ are zero except for the $(1,1)$-entry,
	we observe that
	\begin{align*}
		A E A^\top
		&=
			\begin{pmatrix}
				A_{1,1} & \zerovector_{2m-1}^\top \\
				A_{2:2m,1} & A_{2:2m,2:2m}
			\end{pmatrix}
			\begin{pmatrix}
				0 & \zerovector_{2m-1}^\top \\
				\zerovector_{2m-1} & E_{2:2m,2:2m}
			\end{pmatrix}
			\begin{pmatrix}
				A_{1,1} & A_{2:2m,1}^\top \\
				\zerovector_{2m-1} & A_{2:2m,2:2m}^\top 
			\end{pmatrix} \\
		&=
			\begin{pmatrix}
				0 & \zerovector_{2m-1}^\top \\
				\zerovector_{2m-1} & A_{2:2m,2:2m} E_{2:2m,2:2m} A_{2:2m,2:2m}^\top
			\end{pmatrix}.
	\end{align*}
	Since $E_{2:2m,2:2m}$ is nonnegative-definite and symmetric,
	we have
	\begin{align*}
		\trace(A E A^\top)
		=
			\trace\left(A_{2:2m,2:2m} E_{2:2m,2:2m} A_{2:2m,2:2m}^\top\right)
		\leq
			\| A_{2:2m,2:2m} \|_\op^2 \trace(E).
	\end{align*}
	Thus, \lref{lem89453204920492013} completes the proof.
\end{proof}

\begin{lemma}\label{lem:Variance-of-Yn}
	There exist $C>0$ and $\theta > 0$ such that, for all $k\geq 1$,
	\begin{align*}
		\expect[|Y_{2k+1} - \expect[Y_{2k+1}]|^2]
		\leq
			Ck^{2-\theta}.
	\end{align*}
\end{lemma}
\begin{proof}
	Let $\bar{Y}_n = Y_n - \expect[Y_n]$ for $n\geq 0$ 
	and $\tilde{Z}_{n} = Z_{n} - \expect[Z_{n}|\sigmaField_{n-1}]$ for $n\ge1$.
	Then, by \eqref{eq:Conditional-expectation-of-urnincrement}, we have
	\begin{align*}
		\bar{Y}_{n+1}
		&=
			(Y_n - \expect[Y_n]) + (\expect[Z_{n+1}|\sigmaField_n] - \expect[Z_{n+1}]) + (Z_{n+1} - \expect[Z_{n+1}|\sigmaField_n]) \\
		&= 
			\left(I_{2m} + \frac{H_{n+1}}{n}\right) \bar{Y}_n + \tilde{Z}_{n+1}.
	\end{align*}
	From this recurrence relation, we obtain 
	\begin{align*}
		\expect[\bar{Y}_{n+1} \bar{Y}_{n+1}^\top|\sigmaField_n]
		&=
			\expect\left[\left(\left(I_{2m} + \frac{H_{n+1}}{n}\right) \bar{Y}_n + \tilde{Z}_{n+1}\right)
				\left(\left(I_{2m} + \frac{H_{n+1}}{n}\right) \bar{Y}_n + \tilde{Z}_{n+1}\right)^\top \middle| \sigmaField_n \right] \\
		&=
			\left(I_{2m} + \frac{H_{n+1}}{n}\right) \bar{Y}_n \bar{Y}_n^\top \left(I_{2m} + \frac{H_{n+1}^\top}{n}\right)
			+
			\expect[\tilde{Z}_{n+1}\tilde{Z}_{n+1}^\top|\sigmaField_n].
	\end{align*}
	and taking the expectation, we obtain
	\begin{align*}
		\expect[\bar{Y}_{n+1} \bar{Y}_{n+1}^\top]
		&=
			\left(I_{2m} + \frac{H_{n+1}}{n}\right)
			\expect[\bar{Y}_n \bar{Y}_n^\top]
			\left(I_{2m} + \frac{H_{n+1}^\top}{n}\right)
			+
			\expect[\tilde{Z}_{n+1}\tilde{Z}_{n+1}^\top].
	\end{align*}
	This recurrence relation yields that
	\begin{align*}
		\expect[\bar{Y}_{n+1} \bar{Y}_{n+1}^\top]
		=
			\sum_{i=1}^{n+1}
				\Pi_{n,i}
				\expect[\tilde{Z}_i\tilde{Z}_i^\top]
				\Pi_{n,i}^\top,
	\end{align*}
	which implies
	\begin{align*}
		\expect[| \bar{Y}_{n+1} |^2]
		=
			\trace\expect[\bar{Y}_{n+1} \bar{Y}_{n+1}^\top]
		=
			\sum_{i=1}^{n+1}
				\mathfrak{t}_{n,i},
	\end{align*}
	where
	\begin{align*}
		\mathfrak{t}_{n,i}
		=
			\trace
				\big(
					\Pi_{n,i}
					\expect[\tilde{Z}_i\tilde{Z}_i^\top]
					\Pi_{n,i}^\top
				\big)
		=
			\trace
				\big(
					(P^\top \Pi_{n,i} P) 
					(P^\top \expect[\tilde{Z}_i\tilde{Z}_i^\top] P)
					(P^\top \Pi_{n,i} P)^\top
				\big).
	\end{align*}

	In what follows, we estimate $\expect[| \bar{Y}_{2k+1} |^2]$ with the help of \lref{lem:Pi-trace-inequality}.
	Using $P \mathbf{e}_1 = \frac{1}{\sqrt{2m}}\onevector_{2m}$ and $\tilde{Z}_i^\top \onevector_{2m} = 0$,
	we have
	$
		P^\top \expect[\tilde{Z}_i \tilde{Z}_i^\top] P \mathbf{e}_1
		=
			\zerovector_{2m}
	$.
	Hence, the first column of $P^\top \expect[\tilde{Z}_i \tilde{Z}_i^\top] P$ is the zero vector.
	Similarly, the first row of those matrices is also the zero vector.
	We also see $P^\top \expect[\tilde{Z}_i \tilde{Z}_i^\top] P$ are nonnegative-definite and symmetric matrices.
	Furthermore, we have
	\begin{align*}
		\trace \big(P^\top \expect[\tilde{Z}_i \tilde{Z}_i^\top] P\big)
		=
			\trace \big(\expect[\tilde{Z}_i \tilde{Z}_i^\top] \big)
		=
			\expect[|\tilde{Z}_i|^2]
		\leq
			1.
	\end{align*}
	Applying \lref{lem:Pi-trace-inequality}, for $i=2l-1,2l$ with $k\geq l\geq 1$, we have
	\begin{align*}
		\mathfrak{t}_{2k,i}
		\leq
			C
			\lambda_{k,l}(\rho)^2
            \trace \big(P^\top \expect[\tilde{Z}_i \tilde{Z}_i^\top] P\big)
		\leq
			C
			\lambda_{k,l}(\rho)^2.
	\end{align*}
	In addition, since $\Pi_{2k,2k+1}=I_{2m}$, we have
    $
		\mathfrak{t}_{2k,2k+1}
		=
            \trace \big(\expect[\tilde{Z}_{2k+1} \tilde{Z}_{2k+1}^\top]\big)
		\leq 1
	$.
	Hence, we have
	\begin{align*}
		\expect\left[\left| \bar{Y}_{2k+1} \right|^2\right]
		=
            \sum_{l=1}^{k}
				\mathfrak{t}_{2k,2l-1}
			+
            \sum_{l=1}^{k}
				\mathfrak{t}_{2k,2l}
			+
			\mathfrak{t}_{2k,2k+1}
        \le
			C
            \sum_{l=1}^k
				\lambda_{k,l}(\rho)^2.
	\end{align*}

    From the definition of $\lambda_{k,l}(\rho)$, we see that 
    \begin{align*}
        \sum_{l=1}^{k}
            \lambda_{k,l}(\rho)^2
        &=
            \left(\frac{\Gamma(k+\frac{1}{2}+\rho)}{\Gamma(k+\frac{1}{2})}\right)^2
            \sum_{l=1}^{k}
                \left(\frac{\Gamma(l-\frac{1}{2})}{\Gamma(l-\frac{1}{2}+\rho)}\right)^2.
    \end{align*}
    By \eqref{eq:Gamma-asymptotic}, we have
    \begin{align*}
        \sum_{l=1}^{k}
            \lambda_{k,l}(\rho)^2
        &\leq
			C
			k^{2\rho}
			\sum_{l=1}^{k}
				l^{-2\rho}
		\leq
			C
            \begin{cases}
                k
                &
                \text{if $0 \le \rho < \frac{1}{2}$,} \\
                k \log (k+1)
                &
                \text{if $\rho = \frac{1}{2}$,} \\
                k^{2\rho}
                &
                \text{if $\rho > \frac{1}{2}$.}
            \end{cases}
    \end{align*}
    Since $\rho<1$, there exists $\theta>0$ such that
    $
        \expect[| \bar{Y}_{2k+1} |^2]
        \leq
			Ck^{2-\theta}
	$.
	We complete the proof.
\end{proof}

Finally, we prove \tref{thm:SLLN}.
\begin{proof}[Proof of \tref{thm:SLLN}]
    Let $r= \lfloor \theta^{-1} \rfloor + 2$,
	where $\theta$ is the constant given in \lref{lem:Variance-of-Yn}.
	Then we have $r > \theta^{-1} + 1$.
	Let $\epsilon > 0$ be arbitrary.
	The Chebyshev inequality and \lref{lem:Variance-of-Yn} yield that
	\begin{align*}
		\prob\left(
			\left| \frac{Y_{2k+1}}{2k+1} - \frac{\expect[Y_{2k+1}]}{2k+1} \right| > \epsilon
			\right)
		&\le
			\frac{\expect\left[\left| Y_{2k+1} - \expect[Y_{2k+1}] \right|^2\right]}{\epsilon^2 (2k+1)^2}
		\le
			\frac{C}{\epsilon^2 k^{\theta}}.
	\end{align*}
	Taking subsequences $n = 2k^r+1$ ($k\ge1$) and summing up over $k$, we have
	\begin{align*}
		\sum_{k=1}^{\infty}
			\prob\left(
				\left| \frac{Y_{2k^r+1}}{2k^r+1} - \frac{\expect[Y_{2k^r+1}]}{2k^r+1} \right| > \epsilon
			\right)
		&\le
			\frac{C}{\epsilon^2}
			\sum_{k=1}^{\infty}
				\frac{1}{k^{r \theta}}.
	\end{align*}
	Since $r \theta > 1$, the right-hand side converges.
	Thus, the Borel-Cantelli lemma yields
	\begin{align*}
		\lim_{k\to\infty}
		\left(
			\frac{Y_{2k^r+1}}{2k^r+1}
			-
			\frac{\expect[Y_{2k^r+1}]}{2k^r+1}
		\right)
		=
			0
		\quad \text{a.s.}
	\end{align*}
	Since $k\mapsto 2k^r+1$ is strictly increasing, for any $n\ge3$,
	there exists a unique integer $k_n \ge1$ such that $2k_n^r+1 \le n < 2(k_n + 1)^r+1$.
	Noting that, for $n\ge m$,
	\begin{align*}
		|Y_n - Y_m|
		\le
			\sum_{i=m+1}^{n} |Z_i|
		=
			n - m,
	\end{align*}
	we have
	\begin{align*}
		\left|
			\frac{Y_n}{n}
			-
			\frac{Y_{2k_n^r+1}}{2k_n^r+1}
		\right|
		&=
			\left|
				\frac{Y_n - Y_{2k_n^r+1}}{n}
				+
				Y_{2k_n^r+1} \left(
					\frac{1}{n} - \frac{1}{2k_n^r+1}
				\right)
			\right|\\
		&\le
			\frac{|n - (2k_n^r+1)|}{n} 
			+
			(2k_n^r+1) \frac{|(2k_n^r+1) - n|}{n (2k_n^r+1)} \\
		&\le
			2 \frac{(k_n + 1)^r - k_n^r}{k_n^r}
		=
			2 \left\{
				\left(1 + \frac{1}{k_n}\right)^r  -1
			\right\}.
	\end{align*}
	Thus, we obtain
	\begin{align*}
		\left|
			\frac{Y_n}{n}
			-
			\frac{Y_{2k_n^r+1}}{2k_n^r+1}
		\right|
		\leq
			\frac{C}{k_n}.
	\end{align*}
	As $n\to\infty$, we have $Y_n/n - Y_{2k_n^r+1}/(2k_n^r+1) \to 0$ almost surely.
	We observe that 
	\begin{align*}
		\left|\frac{Y_n}{n} - \frac{\expect[Y_n]}{n}\right|
		&\le
			\left|
				\frac{Y_n}{n} - \frac{Y_{2k_n^r+1}}{2k_n^r+1}
			\right|
			+
			\left|
				\frac{Y_{2k_n^r+1}}{2k_n^r+1} - \frac{\expect[Y_{2k_n^r+1}]}{2k_n^r+1}
			\right|
			+
			\left|
				\frac{\expect[Y_{2k_n^r+1}]}{2k_n^r+1} - \frac{\expect[Y_n]}{n}
			\right|.
	\end{align*} 
	We have already shown that the first and second terms on the right-hand side converge to zero almost surely
	and the third term also converges to zero,
	and then we have $Y_n/n - \expect[Y_n]/n \to 0$ almost surely.
	Combining this with \tref{thm43092043241902494}, we obtain the first assertion.

	The first assertion and $S_n = \mathbf{V} Y_n$
	imply
	\begin{align*}
		\lim_{n\to\infty}
			\frac{S_n}{n}
		=
			\mathbf{V}
			\lim_{n\to\infty}
				\frac{Y_n}{n}
		=
			\frac{1}{2m}
			\mathbf{V}
			\onevector_{2m}
		\quad
		\text{a.s.}
	\end{align*}
	Since $\mathbf{V}\onevector_{2m}=0$, we see the second assertion.
\end{proof}

Furthermore, we can show the following corollary
of the limit of the conditional expectation and covariance of $X_n$.
\begin{corollary}\label{col:Limit-of-conditional-expectation-and-covariance-of-Xn}
	Assume that $\gamma<1$ and $\delta<1$.
	Then, we have
	\begin{align*}
		\lim_{k\to\infty}
			\expect[X_{2k+1} | \sigmaField_{2k}]
		&=
			\bar{v},
		&
		\lim_{k\to\infty}
			\expect[X_{2k} | \sigmaField_{2k-1}]
		&=
			-\bar{v}
		\qquad
		\text{a.s.}
	\end{align*}
	and
	\begin{align*}
		\lim_{n\to\infty}
			\covariance\left(X_n | \sigmaField_{n-1}\right)
		=
			\frac{1}{m} \sum_{i=1}^{m} (v_i - \bar{v})(v_i - \bar{v})^\top 
			\quad \text{a.s.}
	\end{align*}
	In particular, we have
	\begin{align*}
		\lim_{n\to\infty}
			\trace (\covariance(X_n | \sigmaField_{n-1}))
		=
			\frac{1}{m} \sum_{i=1}^{m} |v_i - \bar{v}|^2 >0
			\quad \text{a.s.}
	\end{align*}

	Moreover, the same limits hold without conditioning, i.e.,
	the above convergences remain valid if $\expect[\cdot\mid\sigmaField_{n-1}]$
	and $\covariance(\cdot\mid\sigmaField_{n-1})$ are replaced by
	$\expect[\cdot]$ and $\covariance(\cdot)$, respectively.
\end{corollary}
\begin{proof}
	From $X_n=\mathbf{V} Z_n$ and \eqref{eq:Conditional-expectation-of-urnincrement}, we have
	\begin{align*}
		\expect[X_{n+1} | \sigmaField_n]
		=
			\mathbf{V} \expect[Z_{n+1} | \sigmaField_n]
		=
			\mathbf{V} \frac{H_{n+1} Y_{n}}{n}.
	\end{align*}
	\tref{thm:SLLN} yields
	\begin{align*}
		\lim_{k\to\infty}
			\expect[X_{2k+1} | \sigmaField_{2k}]
		&=
			\mathbf{V}H_\odd
			\frac{1}{2m}
			\onevector_{2m},
		&
		\lim_{k\to\infty}
			\expect[X_{2k} | \sigmaField_{2k-1}]
		&=
			\mathbf{V}H_\even
			\frac{1}{2m}
			\onevector_{2m}
		\quad
		\text{a.s.}
	\end{align*}
	Noting
	\begin{align*}
		\mathbf{V}
		H_\odd
		\onevector_{2m}
		&=
			\mathbf{V}
			\begin{pmatrix}
				2\cdot \onevector_m\\
				\zerovector_{m}
			\end{pmatrix}
		=
			2\sum_{i=1}^{m} v_i
		= 
			2m
			\bar{v},
		&
		\mathbf{V}
		H_\even
		\onevector_{2m}
		&=
			\mathbf{V}
			\begin{pmatrix}
				\zerovector_{m}\\
				2\cdot \onevector_m
			\end{pmatrix}
		=
			2\sum_{i=1}^{m} v_{m+i}
		=
			-
			2m
			\bar{v},
	\end{align*}
	we have the desired result for the conditional expectation.

	Next, we consider the conditional covariance of $X_{n+1}$.
	Since $Z_{n+1} Z_{n+1}^\top = \diag(Z_{n+1})$, we have
	\begin{align*}
		\expect[X_{n+1} X_{n+1}^\top | \sigmaField_n]
		=
			\mathbf{V} \expect[Z_{n+1} Z_{n+1}^\top | \sigmaField_n] \mathbf{V}^\top
		=
			\mathbf{V} \diag(\expect[Z_{n+1} | \sigmaField_n]) \mathbf{V}^\top.
	\end{align*}
	By \eqref{eq:Conditional-expectation-of-urnincrement}, we have
	\begin{align*}
		\covariance(X_{n+1} | \sigmaField_n)
		&=
			\expect[X_{n+1} X_{n+1}^\top | \sigmaField_n]
			-
			\expect[X_{n+1} | \sigmaField_n] \expect[X_{n+1} | \sigmaField_n]^\top \\
		&=
			\mathbf{V} \diag\left(\frac{H_{n+1} Y_{n}}{n}\right) \mathbf{V}^\top
			-
			\mathbf{V} \frac{H_{n+1} Y_{n}}{n} \left(\mathbf{V} \frac{H_{n+1} Y_{n}}{n}\right)^\top.
	\end{align*}
	Thus, by \tref{thm:SLLN}, the following holds almost surely as $n\to\infty$ along odd/even subsequences,
	\begin{align*}
		\covariance(X_n | \sigmaField_{n-1})
		\to
			\begin{cases}
				\frac{1}{m} \sum_{i=1}^{m} v_i v_i^\top -\frac{1}{m^2} \sum_{i=1}^{m} \sum_{j=1}^{m} v_i v_j^\top & \text{if $n$ is odd}, \\
				\frac{1}{m}\sum_{i=1}^{m} v_{m+i} v_{m+i}^\top - \frac{1}{m^2}\sum_{i=1}^{m} \sum_{j=1}^{m} v_{m+i} v_{m+j}^\top & \text{if $n$ is even}.
			\end{cases}
	\end{align*}
	Noting that $v_{m+i} = -v_i$, we have the desired result for the conditional covariance.
	Taking the trace, we obtain the last assertion.

	Finally, the same limits without conditioning can be shown in the same way.
	From the observation for the conditional expectation, we have
	\begin{align*}
		\expect[X_{n+1} X_{n+1}^\top]
		=
			\mathbf{V} \diag\left(\expect[Z_{n+1}]\right) \mathbf{V}^\top.
	\end{align*}
	Thus, we have from \eqref{eq:Conditional-expectation-of-urnincrement} that
	\begin{align*}
		\covariance(X_{n+1})
		&=
			\expect[X_{n+1} X_{n+1}^\top]
			-
			\expect[X_{n+1}] \expect[X_{n+1}]^\top \\
		&=
			\mathbf{V} \diag\left(\frac{H_{n+1} \expect[Y_n]}{n}\right) \mathbf{V}^\top
			-
			\mathbf{V} \frac{H_{n+1} \expect[Y_n]}{n} \left(\mathbf{V} \frac{H_{n+1} \expect[Y_n]}{n}\right)^\top.
	\end{align*}
	By \tref{thm43092043241902494}, we have the desired result.
\end{proof}

\section{Conditional expectation of one-step increments}\label{sec58934850298092}
In this section, we show the next lemma, which plays an important role in the martingale approach.
Recall $T_n$ is defined by \eqref{eq843904820941423}.
\begin{lemma}\label{lem:Conditional-expectation}
	Assume $V_\odd \cap V_\even = \emptyset$.
	For $n\ge1$, we have
	\begin{align*}
		\expect[X_{n+1} | \sigmaField_n]
		&=
			\frac{\delta}{n} S_n 
			+
			\frac{\gamma}{n}(-1)^n T_n
			+
			\left\{
				(1-\gamma) (-1)^n
				-
				\frac{\delta}{n}
				\indicator{\odd}(n)
			\right\}
			\bar{v}.
	\end{align*}
\end{lemma}
\begin{proof}
	From the definition of $X_{n+1}$, we have
	\begin{align*}
		X_{n+1}
		=
			G_{n+1}(X_{\mathcal{U}_n})
		=
			\sum_{k=1}^n G_{n+1}(X_k) \setindicator{\mathcal{U}_n = k}
		=
			\sum_{v \in V}
				\sum_{k=1}^{n}
					G_{n+1}(v)
					\setindicator{X_k = v}
					\setindicator{\mathcal{U}_n = k}.
	\end{align*}
	Hence we have
	\begin{align*}
		\expect[X_{n+1} | \sigmaField_n]
		=
			\sum_{v \in V}
				\sum_{k=1}^{n}
					\expect[G_{n+1}(v)]
					\setindicator{X_k = v}
					\expect[\setindicator{\mathcal{U}_n = k}]
		=
			\frac{1}{n}
			\sum_{v \in V}
				\expect[G_{n+1}(v)]
				Y_n(v).
	\end{align*}
	Recalling the definition of $G_{n+1}(v)$ and using $\indicator{\odd}(n+1) = \frac{1}{2}(1 - (-1)^{n+1})$ and $\indicator{\even}(n+1) = \frac{1}{2}(1 + (-1)^{n+1})$,
	we obtain
	\begin{align*}
		\expect[X_{n+1} | \sigmaField_n]
		&=
			\frac{1}{2n}
			\sum_{v \in V}
				\{
					\expect[G^\odd_{n+1}(v)] \{1 - (-1)^{n+1}\}
					+
					\expect[G^\even_{n+1}(v)] \{1 + (-1)^{n+1}\}
				\}Y_n(v) \\
		&=
			\frac{1}{2n} \{A^+_n + (-1)^n A^-_n\},
	\end{align*}
	where
	\begin{align*}
		A^+_n
		&=
			\sum_{v \in V}
				\left\{
					\expect[G^\odd_{n+1}(v)]
					+
					\expect[G^\even_{n+1}(v)]
				\right\}
				Y_n(v),\\
		A^-_n
		&=
			\sum_{v \in V}
				\left\{
					\expect[G^\odd_{n+1}(v)]
					-
					\expect[G^\even_{n+1}(v)]
				\right\}
				Y_n(v).
	\end{align*}
	Then we see
	\begin{align}
		\label{eq402902394324}
		A^+_n
        &=
            2\delta S_n
            -
            2\delta 
            \indicator{\odd}(n)
			\bar{v},\\
		\label{eq49209423423414}
		A^-_n
        &=
        	2\gamma T_n
			+
			2(1-\gamma)n
			\bar{v}.
	\end{align}
    Hence, we obtain the desired result.

	In what follows, we show \eqref{eq402902394324} and \eqref{eq49209423423414}.
    First, we calculate $A^+_n$.
	From \eqref{eq:Gk-expectation}, we have
	\begin{align*}
		A^+_n
        &=
            \sum_{v \in V_\odd}
                \left\{\alpha v + (1 - \alpha) \bar{v} \right\} Y_n(v)
            +
            \sum_{v \in V_\even}
                \left\{-\beta v + (1 - \beta) \bar{v} \right\} Y_n(v) \\
        &\qquad\qquad
            +
            \sum_{v \in V_\even}
                \left\{\alpha v - (1 - \alpha) \bar{v} \right\} Y_n(v)
            +
            \sum_{v \in V_\odd}
                \left\{-\beta v - (1 - \beta) \bar{v} \right\} Y_n(v) \\
        &=
            \alpha \left(\sum_{v\in V_\odd} + \sum_{v \in V_\even}\right) 
                v Y_n(v)
            -
            \beta \left(\sum_{v \in V_\even} + \sum_{v \in V_\odd}\right)
                v Y_n(v) \\
        &\qquad\qquad
            +
            (1 - \alpha) \bar{v} \left(\sum_{v \in V_\odd} - \sum_{v \in  V_\even}\right)
                Y_n(v)
            +
            (1 - \beta) \bar{v} \left(\sum_{v \in V_\even} - \sum_{v \in V_\odd}\right)
                Y_n(v).
	\end{align*}
	Hence we see \eqref{eq402902394324} by applying \eqref{eq453849023809481} and
	\begin{align*}
		\left(\sum_{v \in V_\odd} - \sum_{v \in  V_\even}\right) 
            Y_n(v)
        &=
            \indicator{\odd}(n).
	\end{align*}
	We can arrive at \eqref{eq49209423423414} in the same way as $A^+_n$
	by using \eqref{eq843904820941423} and
	\begin{align*}
		\left(\sum_{v \in V_\odd} + \sum_{v \in  V_\even}\right) 
            Y_n(v)
		&=
			n.
	\end{align*}
	The proof is complete.
\end{proof}

\section{Centered second moment}\label{sec:Centered-second-moment}
In this section, we show \tref{thm:Asymptotic-second-moment}.
To this end, we need to investigate the asymptotic behavior of the second moments of
$S_n$ and $T_n$. We write the centered variables as
\begin{align*}
	\bar{S}_n = S_n - \expect[S_n],
	\qquad
	\bar{T}_n = T_n - \expect[T_n],
	\qquad
	\bar{X}_n = X_n - \expect[X_n].
\end{align*}
The constant $\sigma^2$ given by \eqref{eq489208401} is well-defined
due to \corref{col:Limit-of-conditional-expectation-and-covariance-of-Xn}.
Therefore, we reduce \tref{thm:Asymptotic-second-moment} to the following theorem.
\begin{theorem}\label{thm:Asymptotic-second-moment-for-Sn-Tn}
	Assume that $\gamma<1$ and $\delta<1$. Then, we have
	\begin{align*} 
		&\expect[|\bar{S}_n|^2]
		\sim
			\begin{cases}
				\frac{\sigma^2}{1-2\delta}n & \text{if $\delta<\frac{1}{2}$}, \\
				\sigma^2 \, n\log n & \text{if $\delta=\frac{1}{2}$}, \\
				C_{\alpha, \beta} n^{2\delta}  & \text{if $\delta>\frac{1}{2}$},
			\end{cases}
		&
		&\expect[|\bar{T}_n|^2]
		\sim
			\begin{cases}
				\frac{\sigma^2}{1-2\gamma}n & \text{if $\gamma<\frac{1}{2}$}, \\
				\sigma^2 \, n\log n & \text{if $\gamma=\frac{1}{2}$}, \\
				C^\prime_{\alpha, \beta} n^{2\gamma}  & \text{if $\gamma>\frac{1}{2}$},
			\end{cases}
	\end{align*}
	where $C_{\alpha, \beta}$ and $C^\prime_{\alpha, \beta}$
	are positive constants depending on $\alpha$ and $\beta$.
\end{theorem}

Let $s_n=\expect[|\bar{S}_n|^2]$, $t_n=\expect[|\bar{T}_n|^2]$, 
$u_n=\expect[\bar{S}_n^\top \bar{T}_n]$, $\sigma^2_n = \expect[|\bar{X}_n|^2]$.
From the definitions, we observe that 
\begin{align*}
	|S_n|, |T_n|
	\le
		\sum_{k=1}^{n} |X_k|
	\le
		n \cdot \sup_{v\in V} |v|,
\end{align*}
which implies that there exists a constant $C>0$ such that, for $n\ge1$,
\begin{align}\label{eq489329084109}
	|s_n|
	&\leq
		C n^2, 
	&
	|t_n|
	&\leq
		C n^2,
	&
	|u_n|
	&\leq
		C n^2.
\end{align}
We also see
\begin{align}
	\label{eq48390284109824}
	\sigma^2_n
	=
		\expect[|X_n|^2]
		-
		|\expect[X_n]|^2
	\le
		\expect[|X_n|^2]
	\le
		\sup_{v\in V} |v|^2.
\end{align}

Through a sequence of lemmas, we establish \tref{thm:Asymptotic-second-moment-for-Sn-Tn}.
\begin{lemma}\label{lem:Recurrence-of-second-moment}
	For $n\ge 1$, we have
	\begin{align*}
		s_{n+1}
		&=
			\left(1+\frac{2\delta}{n}\right)s_n
			+ 
			\frac{2\gamma}{n}(-1)^n u_n 
			+
			\sigma^2_{n+1}, \\
		t_{n+1}
		&=
			\left(1+\frac{2\gamma}{n}\right)t_n
			+ 
			\frac{2\delta}{n}(-1)^n u_n 
			+
			\sigma^2_{n+1}, \\
		u_{n+1}
		&=
			\left(1+\frac{\delta+\gamma}{n}\right) u_n 
			+ 
			\frac{(-1)^n}{n}\left(\delta s_n+ \gamma t_n\right)
			+
			(-1)^n \sigma^2_{n+1}.
	\end{align*}
\end{lemma}
\begin{proof}
	From \lref{lem:Conditional-expectation},
	the conditional expectation of $\bar{X}_{n+1}$ is rewritten as
	\begin{align*}
		\expect[\bar{X}_{n+1}|\sigmaField_n]
		=
		\expect[X_{n+1}|\sigmaField_n] - \expect[X_{n+1}]
		=
			\frac{\delta}{n} \bar{S}_n + \frac{\gamma}{n} (-1)^n \bar{T}_n.
	\end{align*}
	By the definitions of $\bar{S}_n$ and $\bar{T}_n$, we have
	\begin{align*}
		\bar{S}_{n+1} = \bar{S}_n + \bar{X}_{n+1},
		\quad
		\bar{T}_{n+1} = \bar{T}_n + (-1)^n \bar{X}_{n+1}
	\end{align*}
	and
	\begin{align*}
		\bar{S}_{n+1}^\top \bar{T}_{n+1}
		=
			\bar{S}_n^\top \bar{T}_n 
			+
			(-1)^n \bar{S}_n^\top \bar{X}_{n+1} 
			+ 
			\bar{X}_{n+1}^\top \bar{T}_n 
			+ 
			(-1)^n |\bar{X}_{n+1}|^2.
	\end{align*}
	Taking expectations of $|\bar{S}_{n+1}|^2, |\bar{T}_{n+1}|^2$,
	and $\bar{S}_{n+1}^\top \bar{T}_{n+1}$ yields
	\begin{align*}
		s_{n+1}
		&=
			s_n 
			+
			2 \expect[\bar{S}_n^\top \bar{X}_{n+1}] 
			+ 
			\sigma^2_{n+1}, \\
		t_{n+1}
		&=
			t_n 
			+
			2 (-1)^n \expect[\bar{T}_n^\top \bar{X}_{n+1}] 
			+ 
			\sigma^2_{n+1}, \\
		u_{n+1}
		&=
			u_n 
			+
			(-1)^n \expect[\bar{S}_n^\top \bar{X}_{n+1}] 
			+ 
			\expect[\bar{T}_n^\top \bar{X}_{n+1}] 
			+ 
			(-1)^n \sigma^2_{n+1}.
	\end{align*}
	By the conditional expectation of $\bar{X}_{n+1}$, we have
	\begin{align*}
		\expect[\bar{S}_n^\top \bar{X}_{n+1}]
		&=
			\frac{\delta}{n} s_n
			+ 
			\frac{\gamma}{n} (-1)^n u_n,
		&
		\expect[\bar{T}_n^\top \bar{X}_{n+1}]
		&=
			\frac{\delta}{n}  u_n
			+ 
			\frac{\gamma}{n} (-1)^n t_n.
	\end{align*}
	Substituting these results into the recurrence relations, we obtain the desired recursions.
\end{proof}

We study $s_n$ and $t_n$ using the recurrence relation in \lref{lem:Recurrence-of-second-moment}.
Given their similarity, we restrict our attention mainly to $s_n$.
To this end, set $\eta_3=1$ and
\begin{align*}
	\eta_{n+1}
	=
		\left(1+\frac{2\delta}{n}\right) \eta_n
	\qquad (n\ge3).
\end{align*}
Since $1+\frac{2\delta}{n}\geq 1+\frac{2(-1)}{3}>0$ for $n\geq 3$,
we have $\eta_n>0$.
It follows from  \eqref{eq:Gamma-asymptotic} that
\begin{align}\label{eq:etan-asymptotic}
	\eta_n
	=
		\left\{
			\prod_{k=3}^{n-1} \left(1+\frac{2\delta}{k}\right)
		\right\}
		\eta_3
	=
		\frac{\Gamma(n+2\delta)}{\Gamma(3+2\delta)}
		\frac{\Gamma(3)}{\Gamma(n)}
	\sim
		\frac{2}{\Gamma(3+2\delta)} n^{2\delta}.
\end{align}
Dividing the recurrence relation for $s_n$ by $\eta_{n+1}$ yields
\begin{align*}
	\frac{s_{n+1}}{\eta_{n+1}} - \frac{s_n}{\eta_n}
	=
		\frac{1}{\eta_{n+1}}
			\left(\frac{2\gamma}{n}(-1)^n u_n + \sigma^2_{n+1}\right).
\end{align*}
For $n\geq 4$, by summing this from $n=3$ to $n-1$, we see
\begin{align}\label{eq4832982}
	\frac{s_n}{\eta_n}
	=
		\frac{s_3}{\eta_3}
		+
		\sum_{k=3}^{n-1}
			\left(
				\frac{s_{k+1}}{\eta_{k+1}} - \frac{s_k}{\eta_k}
			\right)
	=
		s_3
		+
		\sum_{k=3}^{n-1}
			\frac{1}{\eta_{k+1}}
			\left(\frac{2\gamma}{k}(-1)^k u_k + \sigma^2_{k+1}\right).
\end{align}
Hence for $n\geq 4$, $s_n$ can be rewritten as
\begin{align}\label{eq:s_n-Explicity-Solution}
	s_n 
	= 
			s_3\eta_n
			+
			\eta_n\sum_{k=4}^{n} \frac{\sigma^2_{k}}{\eta_k} 
			+ 
			2\gamma
			\eta_n
			\Sigma_n,
\end{align}
where 
\begin{align*}
	\Sigma_n 
	=
		\sum_{k=3}^{n-1} \frac{(-1)^k}{k} \frac{u_k}{\eta_{k+1}}.
\end{align*}

To analyze $s_n$, we first estimate $\Sigma_n$.
\begin{lemma}\label{lem:Inequality-for-second-term-in-s_n}
	There exists a constant $C>0$ such that, for $n\ge4$,
	\begin{align}
		\label{eq4328990482309}
		|\Sigma_n|
		\le
			C
			\left(
				\frac{|u_n|}{n\eta_{n+1}}
				+
				\sum_{k=3}^{n-1}
					\frac{|u_k| + |s_k| + |t_k| + k}{k^2\eta_{k+2}}
			\right).
	\end{align}
\end{lemma}
\begin{proof}
	Let  $A_k=\frac{u_k}{k\eta_{k+1}}$ for $k\geq 2$ and 
	$B_2=0$ and $B_n = \sum_{k=3}^{n} (-1)^k$ for $n\geq 3$.
	By the summation by parts, we have
	\begin{align*}
		\Sigma_n
		=
			\sum_{k=3}^{n-1}
				A_k \{B_k - B_{k-1}\}
		=
			A_nB_{n-1}
			-
			\sum_{k=3}^{n-1}
				\{A_{k+1}-A_k\}
				B_k.
	\end{align*}
	Note that
	$
		|A_nB_{n-1}|
		\leq
			|A_n|
	$.
	From the definition of $\eta_n$, we have
	\begin{align*}
		A_{k+1}-A_k
		&=
			\frac{u_{k+1}}{(k+1)\eta_{k+2}} - \frac{u_k}{k\eta_{k+2}}\left(1+\frac{2\delta}{k+1}\right) \\
		&=
			\frac{k(u_{k+1} - u_k)-(1+2\delta)u_k}{k(k+1)\eta_{k+2}}.
	\end{align*}
	From the recurrence relation of $u_n$ in \lref{lem:Recurrence-of-second-moment}, we have
	\begin{align*}
		k(u_{k+1} - u_k) - (1+2\delta)u_k
		&=
			(\delta + \gamma) u_k + (-1)^k (\delta s_k + \gamma t_k) + (-1)^k k \sigma^2_{k+1}
			-
			(1+2\delta) u_k \\
		&=
			-(1+\delta-\gamma) u_k + (-1)^{k} (\delta s_k + \gamma t_k) + (-1)^{k} k \sigma^2_{k+1}.
	\end{align*}
	Using \eqref{eq489329084109} and \eqref{eq48390284109824}, we have
	\begin{align*}
		\left|
			\sum_{k=3}^{n-1}
				\{A_{k+1}-A_k\}
				B_k
		\right|
		\leq
			C
			\sum_{k=3}^{n-1}
				\frac{|u_k| + |s_k| + |t_k| + k}{k^2 \eta_{k+2}}.
	\end{align*}
	Here $C>0$ is a constant.
	Hence we obtain the desired result.
\end{proof}

The following lemma follows directly from \lref{lem:Inequality-for-second-term-in-s_n}.
\begin{lemma}\label{lem:Rough-upper-bound-for-second-moment}
	We have
	\begin{align*}
			&s_n
			=
				\begin{cases}
					O(n^{\max\{1, 2\delta\}}) & \text{if $\delta\neq 1/2$}, \\
					O(n\log n) & \text{if $\delta= 1/2$},
				\end{cases} 
			&
			t_n
			=
				\begin{cases}
					O(n^{\max\{1, 2\gamma\}}) & \text{if $\gamma\neq 1/2$},\\
					O(n\log n) & \text{if $\gamma= 1/2$}.
				\end{cases}
		\end{align*}
\end{lemma}
\begin{proof}
	It follows from \lref{lem:Inequality-for-second-term-in-s_n} and \eqref{eq489329084109} that
	\begin{align*}
		|\Sigma_n|
		\leq
			C
			\left(
				\frac{|u_n|}{n\eta_{n+1}}
				+
				\sum_{k=3}^{n-1}
					\frac{|u_k| + |s_k| + |t_k| + k}{k^2\eta_{k+2}}
			\right)
		\leq
			C
			\left(
				\frac{n}{\eta_{n+1}}
				+
				\sum_{k=3}^{n-1}
					\frac{1}{\eta_{k+2}}
			\right),
	\end{align*}
	which implies
	\begin{align*}
		|s_n|
		\le
			\eta_n |s_3|
			+
			\eta_n \sum_{k=4}^{n} \frac{C}{\eta_k}
			+
			C\eta_n|\Sigma_n| 
		\le
			C \eta_n
			\left(
				\sum_{k=3}^{n+1} \frac{1}{\eta_k}
				+
				\frac{n}{\eta_{n+1}}
			\right).
	\end{align*}
	For the case $\delta \le 1/2$, \eqref{eq:etan-asymptotic} implies
	\begin{align*}
		\sum_{k=3}^{n+1}
			\frac{1}{\eta_k}
		\sim 
			\begin{cases}
				\frac{\Gamma(3+2\delta)}{2} \frac{n^{1-2\delta}}{1-2\delta} & \text{if $\delta<1/2$}, \\
				3 \log n & \text{if $\delta=1/2$}, \\
			\end{cases}  
	\end{align*}
	and, for the case $\delta>1/2$, the summation converges.
	Combining this and \eqref{eq:etan-asymptotic}, we obtain 
	\begin{align*}
		|s_n|
		\leq
			C
			\eta_n 
			\left(
				\sum_{k=3}^{n+1} \frac{1}{\eta_k}
				+
				\frac{n}{\eta_{n+1}}
			\right)
		\leq
			C
			\begin{cases}
				n^{2\delta} & \text{if $\delta>1/2$}, \\
				n\log n & \text{if $\delta=1/2$}, \\
				n & \text{if $\delta<1/2$}.
			\end{cases}
	\end{align*}
	Thus, we have the bound for $s_n$.

	By the same argument, with $\delta$ and $s_n$ replaced by $\gamma$ and $t_n$,
	respectively, we obtain the bound for $t_n$.
\end{proof}

\begin{proof}[Proof of \tref{thm:Asymptotic-second-moment-for-Sn-Tn}]
	Since we can give the asymptotic behavior of
	$s_n=\expect[|\bar{S}_n|^2]$ for $\delta \le 1/2$
	and $t_n=\expect[|\bar{T}_n|^2]$ for $\gamma \le 1/2$ in the same way,
	we consider $s_n$ for $\delta \le 1/2$ only.
	Recall the explicit expression \eqref{eq:s_n-Explicity-Solution} of $s_n$
	and estimate each term.
	By \eqref{eq:etan-asymptotic}, the first term of \eqref{eq:s_n-Explicity-Solution}
	is $o(n)$ when $\delta<1/2$ and $o(n\log n)$ when $\delta=1/2$.
	Using \eqref{eq489208401} and \eqref{eq:etan-asymptotic},
	the second term of \eqref{eq:s_n-Explicity-Solution} is estimated as
	\begin{align*}
		\eta_n\sum_{k=4}^{n} \frac{\sigma^2_{k}}{\eta_k} 
		\sim
			\begin{cases}
				\frac{\sigma^2}{1-2\delta}n & \text{if $\delta<\frac{1}{2}$}, \\
				\sigma^2 \, n\log n & \text{if $\delta=\frac{1}{2}$}.
			\end{cases}
	\end{align*}
	Hence, it suffices to show that the last term of \eqref{eq:s_n-Explicity-Solution} is estimated as
	\begin{align*}
		\eta_n
		|\Sigma_n|
		=
			\begin{cases}
				o(n) & \text{if $\delta<\frac{1}{2}$},\\
				o(n\log n) & \text{if $\delta=\frac{1}{2}$}.
			\end{cases}
	\end{align*}
	
	We consider the case $\delta<1/2$
	and estimate the right-hand side of \eqref{eq4328990482309}.
	From \lref{lem:Rough-upper-bound-for-second-moment}, we have $s_n=O(n)$.
	Picking $\rho\in (\max\{\gamma, 1/2, 1/2+\delta\}, 1)$, we have $t_n=O(n^{2\rho})$ from \lref{lem:Rough-upper-bound-for-second-moment}.
	By the Cauchy-Schwarz inequality, we obtain
	$
		|u_n|
		\le
			\sqrt{s_n t_n}
		\leq
			Cn^{\rho+\frac{1}{2}}
	$.
	These results yield that $s_n, t_n, |u_n| \leq Cn^{2\rho}$.
	From this, we have
	\begin{align*}
		\eta_n
		\frac{|u_n|}{n\eta_{n+1}} 
		\leq 
			C
			n^{2\rho-1}
		=
			o(n),
	\end{align*}
	which corresponds to the first term in the right-hand side of \eqref{eq4328990482309}.
	Next we estimate the second term in the right-hand side of \eqref{eq4328990482309}.
	Using the estimates above and \eqref{eq:etan-asymptotic}, 
	the summand in \eqref{eq4328990482309} is bounded by
	\begin{align*}
			\frac{|u_k| + |s_k| + |t_k| +k}{k^2\eta_{k+2}}
			\leq
				Ck^{2\rho}k^{-2-2\delta}
			=
				Ck^{2\rho-2\delta-2}.
	\end{align*}
	Thus, by \eqref{eq:etan-asymptotic} and the above bound, we have
	\begin{align*}
		\eta_n\sum_{k=3}^{n-1}
				\frac{|u_k| + |s_k| + |t_k| +k}{k^2\eta_{k+2}}
		\leq
			Cn^{2\delta} \cdot n^{2\rho-2\delta-1}
		=
			Cn^{2\rho-1}
		=
			o(n).
	\end{align*}
	These results yield the desired result.

	Next, we consider the case $\delta=1/2$.
	We see $\frac{1}{2}-\frac{1}{m-1}\leq \gamma\leq \frac{1}{2}$ due to \eqref{eq4539028429084}.
	From \lref{lem:Rough-upper-bound-for-second-moment}, 
	we have $s_n=O(n\log n)$ and $t_n=O(n\log n)$.
	The Cauchy-Schwarz inequality yields
	$
		|u_n|
		\le
			\sqrt{s_n t_n}
		\leq
			Cn\log n
	$.
	Using these results and \eqref{eq:etan-asymptotic},
	the summand on the right-hand side of \eqref{eq4328990482309} is bounded by
	\begin{align*}
		\frac{|u_k| + |s_k| + |t_k| +k}{k^2\eta_{k+2}}
		\leq
			C
			(k\log k)(k^{-2-1})
		=
			C
			k^{-2}\log k.
	\end{align*}
	Hence the sum on the right-hand side of \eqref{eq4328990482309} is bounded uniformly in $n$.
	Therefore, we have
	\begin{align*}
		\eta_n
		|\Sigma_n|
		\leq
			C
			n
			\left\{
				\frac{n\log n}{n\cdot n}
				+
				1
			\right\}
		\leq
			C(\log n+n).
	\end{align*}
	These results yield the desired result.

	Since we can deal with $s_n$ for $\delta > 1/2$ and $t_n$ for $\gamma > 1/2$ in the same way,
	we consider $s_n$ for $\delta > 1/2$.
	From \eqref{eq4539028429084}, $\delta > 1/2$ implies $\gamma<1/2$.
	In this case, we show that
	\begin{align*}
		\lim_{n\to\infty}
			\frac{s_n}{n^{2\delta}}
		=
			C_{\alpha,\beta}>0.
	\end{align*}
	We see $s_n = O(n^{2\delta})$, $t_n = O(n)$ and $u_n = O(n^{\delta+\frac12})$
	from \lref{lem:Rough-upper-bound-for-second-moment} and the Cauchy-Schwarz inequality.
	We can observe from \eqref{eq4832982} that
	\begin{align*}
		s_n 
		=
			\eta_n
			\left\{
				s_3
				+
				\sum_{k=3}^{n-1}
					\frac{1}{\eta_{k+1}}
					\left(\frac{2\gamma}{k}(-1)^k u_k + \sigma^2_{k+1}\right)
			\right\}.
	\end{align*}
	From $u_n = O(n^{\delta+\frac12})$ and \eqref{eq:etan-asymptotic}, the summand is bounded by
	\begin{align*}
		\left|
			\frac{1}{\eta_{k+1}}
			\left(\frac{2\gamma}{k}(-1)^k u_k + \sigma^2_{k+1}\right)
		\right|
		\leq
			C
			k^{-2\delta}
			\left\{
				k^{\left(\delta+\frac12\right)-1}+1
			\right\}
		\leq
			Ck^{-\delta-\frac12}.
	\end{align*}
	Since $\delta>1/2$, the summation converges.
	Therefore, it follows from \eqref{eq:etan-asymptotic} that
	\begin{align*}
		\lim_{n\to\infty}
			\frac{s_n}{n^{2\delta}}
		=
			C_{\alpha,\beta},
			\quad
			\text{where}
			\quad
			C_{\alpha,\beta}
			=
				\frac{\Gamma(3)}{\Gamma(3+2\delta)} 
					\left\{
						s_3
						+
						\sum_{k=3}^{\infty} \frac{1}{\eta_{k+1}} 
						\left(\frac{2\gamma}{k}(-1)^k u_k + \sigma^2_{k+1}\right)
					\right\}.
	\end{align*}

	Next, we claim that $C_{\alpha,\beta}>0$.
	Setting $s^{\prime}_n = s_n - \frac{\gamma}{\delta} t_n$, and using the recurrence relations of $s_n$ and $t_n$, we obtain
	\begin{align*}
		s^{\prime}_{n+1}
		=
			\left(1+\frac{2\delta}{n}\right) s^{\prime}_n 
			+ 
			\left(1-\frac{\gamma}{\delta}\right)
			\left(\sigma_{n+1}^2+\frac{2\gamma}{n}t_n\right).
	\end{align*}

	Here we see that $s^{\prime}_n\ge0$ for all $n\ge3$ as follows.
	When $\gamma<0$, $s^{\prime}_n = s_n +\frac{|\gamma|}{\delta}t_n \ge0$.
	We can rewrite this recurrence relation as
	\begin{align*}
		\frac{s^{\prime}_{n+1}}{\eta_{n+1}} - \frac{s^{\prime}_{n}}{\eta_{n}}
		=
			\frac{1}{\eta_{n+1}}
			\left(1-\frac{\gamma}{\delta}\right)
			\left(\sigma_{n+1}^2 + \frac{2\gamma}{n} t_n\right).
	\end{align*}
	This implies that $s^{\prime}_n/\eta_n$ is monotonically increasing when $0\le\gamma<1/2$.
	We observe from the recurrence relations of $s_n$ and $t_n$ that
	\begin{align*}
		s_3
		&=
			(1+\delta)(1+2\delta)s_1 
			-
			2\gamma(1+\delta) u_1
			+
			\gamma u_2 
			+
			(1+\delta)\sigma^2_2
			+
			\sigma^2_3, \\
		t_3
		&=
			(1+\gamma)(1+2\gamma) t_1
			-
			2\delta(1+\gamma) u_1
			+
			\delta u_2
			+
			(1+\gamma)\sigma^2_2
			+
			\sigma^2_3
	\end{align*}
	and thus
	\begin{align*}
		s^{\prime}_3
		=
			s_3 - \frac{\gamma}{\delta} t_3 
		&=
			(1+\delta)(1+2\delta)s_1 - \frac{\gamma}{\delta} (1+\gamma)(1+2\gamma) t_1 
			+
			2\gamma(\gamma-\delta) u_1 \\
		&\phantom{=}
			\qquad
			+
			\left\{(1+\delta) - \frac{\gamma}{\delta}(1+\gamma)\right\}\sigma^2_2
			+
			\left(1 - \frac{\gamma}{\delta}\right) \sigma^2_3.
	\end{align*}
	Noting that $\bar{S}_1 = \bar{X}_1$ and $\bar{T}_1 = \bar{X}_1$, we have $s_1 = t_1 = \sigma^2_1 >0$ and $u_1 = \sigma^2_1>0$,
	which imply
	\begin{align*}
		s^{\prime}_3
		&=
			\left\{
				(1+\delta)(1+2\delta) - \frac{\gamma}{\delta} (1+\gamma)(1+2\gamma) 
				+
				2\gamma(\gamma-\delta) 
			\right\}\sigma^2_1 \\
		&\phantom{=}
			\qquad
			+
			\left\{(1+\delta) - \frac{\gamma}{\delta}(1+\gamma)\right\}\sigma^2_2
			+
			\left(1 - \frac{\gamma}{\delta}\right) \sigma^2_3 \\
		&=
			\frac{\delta-\gamma}{\delta}
			\left\{
				(1 + 3(\delta+\gamma) + 2(\delta^2 + \gamma^2))\sigma^2_1
				+
				(1+\delta+\gamma)\sigma^2_2
				+
				\sigma^2_3
			\right\}.
	\end{align*}
	Since $\delta>1/2$, $0\le \gamma <1/2$ and $\sigma^2_n\ge0$ for all $n$, $s^{\prime}_3/\eta_3$ is positive.
	Thus,
	we obtain that $s^{\prime}_n/\eta_n$ is positive for $n\ge3$ since this is monotonically increasing.

	Recall $\delta > 1/2$ implies $\gamma<1/2$.
	Then, it follows from \eqref{eq489208401} and $t_n \sim \sigma^2 n/(1-2\gamma)$ 
	for $\gamma<1/2$,
	which has already been shown, that 
	\begin{align*}
		\lim_{n\to\infty}
			\left(\sigma_{n+1}^2 + \frac{2\gamma}{n} t_n\right)
		=
			\sigma^2 + \frac{2\gamma\sigma^2}{1-2\gamma}
		=
			\frac{\sigma^2}{1-2\gamma}
		>
			0.
	\end{align*}
	Therefore, there exists $N$ so that for all $n\ge N$,
	\begin{align*}
		\sigma_{n+1}^2 + \frac{2\gamma}{n} t_n
		\geq
			\frac{\sigma^2}{2(1-2\gamma)}
		>
			0.
	\end{align*}
	By the recurrence relation of $s^{\prime}_n$ and the nonnegativity of $s^{\prime}_N/\eta_N$, we obtain for $n>N$,
	\begin{align*}
		\frac{s^{\prime}_n}{\eta_n}
		=
			\frac{s^{\prime}_N}{\eta_N} 
			+ 
			\left(1-\frac{\gamma}{\delta}\right)
			\sum_{k=N}^{n-1}
				\frac{1}{\eta_{k+1}} 
				\left(\sigma_{k+1}^2 + \frac{2\gamma}{k} t_k\right) 
		\geq
			\left(1-\frac{\gamma}{\delta}\right)
			\frac{\sigma^2}{2(1-2\gamma)}
			\sum_{k=N}^{n-1}
				\frac{1}{\eta_{k+1}}.
	\end{align*}
	Taking the limit inferior on both sides and using
	\eqref{eq:etan-asymptotic} yield
	\begin{align*}
		\liminf_{n\to\infty} \frac{s^{\prime}_n}{n^{2\delta}}
		\geq
			\frac{\Gamma(3)}{\Gamma(3+2\delta)}
			\left(1-\frac{\gamma}{\delta}\right)
			\frac{\sigma^2}{2(1-2\gamma)}
			\sum_{k=N}^{\infty}
				\frac{1}{\eta_{k+1}}
		>0.
	\end{align*}
	As we have already shown that $s_n/n^{2\delta}$ converges to $C_{\alpha,\beta}$ and $t_n = o(n^{2\delta})$, 
	it follows from $\liminf_{n\to\infty} s^{\prime}_n/n^{2\delta}>0$ that
	\begin{align*}
		C_{\alpha,\beta}
		=
			\lim_{n\to\infty} \frac{s_n}{n^{2\delta}}
		=
			\lim_{n\to\infty} \left(\frac{s_n}{n^{2\delta}} - \frac{\gamma}{\delta}\frac{t_n}{n^{2\delta}} \right)
		=
			\liminf_{n\to\infty} \frac{s^{\prime}_n}{n^{2\delta}}
		>
			0.
	\end{align*}
	Hence, we complete the proof of $s_n$ for the case $\delta>1/2$.
\end{proof}

\section{Some limit theorems with the martingale approach}\label{sec:limit-theorems}
\subsection{Construction of martingale}
In order to investigate the asymptotic behavior of $S_n$, 
we find a $d$-dimensional stochastic process $\{M_n\}_{n\geq 0}$ 
of the form $M_n=\zeta_n^{-1} S_n - R_n$,
where $\{\zeta_n\}_{n\geq 0}$ is a deterministic positive sequence
and $\{R_n\}_{n\geq 0}$ is a predictable process with respect to $\{\sigmaField_n\}_{n\geq 0}$.
Here, $\{\sigmaField_n\}_{n\geq 0}$ is the filtration introduced in \secref{sec9092423242092}.
Under this setting, we find $\{\zeta_n\}_{n\geq 0}$ and $\{R_n\}_{n\geq 0}$
so that $\{M_n\}_{n\geq 0}$ is a martingale with respect to $\{\sigmaField_n\}_{n\geq 0}$
and that $R_n$ does not contain the term $S_n$.
Since $\{M_n\}_{n\geq 0}$ should be a martingale, we expect
\begin{align*}
	0
	=
		\expect[M_{n+1}|\sigmaField_n]-M_n
	=
		\zeta_{n+1}^{-1} \{S_n+\expect[X_{n+1}|\sigmaField_n]\}-R_{n+1}
		-\zeta_n^{-1}S_n+R_n,
\end{align*}
which implies
\begin{align*}
	R_{n+1}-R_n
	=
		\frac{1}{\zeta_{n+1}}
		\left\{
			\expect[X_{n+1}|\sigmaField_n]
			+
			\left(
				1-\frac{\zeta_{n+1}}{\zeta_n}
			\right)
			S_n
		\right\}.
\end{align*}
Comparing this expression with \lref{lem:Conditional-expectation},
if $1-\frac{\zeta_{n+1}}{\zeta_n}=-\frac{\delta}{n}$ holds,
then $R_n$ does not contain the term $S_n$.

Based on the above considerations, we define $\{\zeta_n\}_{n\geq 2}$ and $\{R_n\}_{n\geq 2}$ by
\begin{align*}
	\zeta_2
	&=
		1,
	&
	\zeta_{n+1}
	&=
		\left(1 + \frac{\delta}{n}\right)
		\zeta_n
	\qquad
	(n\geq 2),\\
	R_2
	&=
		0,
	&
	R_{n+1}
	&=
		R_n
		+
		\frac{1}{\zeta_{n+1}}
		\left(\expect[X_{n+1}|\sigmaField_n] - \frac{\delta}{n} S_n\right)
	\qquad
	(n\geq 2).
\end{align*}
Using $\{\zeta_n\}_{n\geq 2}$ and $\{R_n\}_{n\geq 2}$, 
we define a $d$-dimensional stochastic process $\{M_n\}_{n\ge2}$ by
\begin{align*}
	M_n 
	= 
		\zeta_n^{-1} S_n - R_n - S_2.
\end{align*}
We now define all sequences only for $n\geq 2$
in order to treat all cases $-1\leq \delta<1$ simultaneously.
Note that $1 + \frac{\delta}{1}=0$ when $\delta=-1$.
Note that the term $S_2$ is included so that $M_2=0$.

We now examine properties of $\{M_n\}_{n\ge2}$.
From the definition above, $\{M_n\}_{n\ge2}$ is 
a square-integrable stochastic process with $M_2=0$.
In addition, we see it is a martingale with respect to the filtration $\{\sigmaField_n\}_{n\ge2}$
because $S_2$ is $\sigmaField_2$-measurable and, for $n\geq 2$,
\begin{align*}
	\expect[\zeta_{n+1}^{-1} S_{n+1} - R_{n+1}|\sigmaField_n] 
	&=
		\zeta_{n+1}^{-1}\left\{S_n + \expect[X_{n+1}|\sigmaField_n]\right\}\\
	&\phantom{=}
		\qquad
		\qquad
		-
		\left\{
			R_n
			+
			\frac{1}{\zeta_{n+1}}
			\left(\expect[X_{n+1} | \sigmaField_n] - \frac{\delta}{n} S_n\right)
		\right\}\\
	&=
		\frac{1}{\zeta_{n+1}}
		\left(1 + \frac{\delta}{n}\right)
		S_n 
		-
		R_n\\
	&=
		\zeta_n^{-1} S_n - R_n.
\end{align*}
To analyze the asymptotic behavior of $M_n$, we consider its predictable quadratic variation.
It is given by $\langle M \rangle_2=0$ and, for $n\geq 3$,
\begin{align}\label{eq:Predictable-quad-variation-of-M_n}
	\langle M \rangle_n 
	= 
		\sum_{k=3}^n \expect[\Delta M_k (\Delta M_k)^\top |\sigmaField_{k-1}] 
	=
		\sum_{k=3}^n \zeta_k^{-2} \covariance(X_k|\sigmaField_{k-1}).
\end{align}
Here we used
\begin{gather}
	\label{eq840928409231}
	\begin{aligned}
		\Delta M_k
		&=
			M_k-M_{k-1}\\
		&=
			\zeta_k^{-1} S_k - \zeta_k^{-1} \left(1 + \frac{\delta}{k-1}\right)S_{k-1} - (R_k - R_{k-1}) \\
		&=
			\zeta_k^{-1} \left(X_k - \expect[X_k|\sigmaField_{k-1}]\right).
	\end{aligned}
\end{gather}

\subsection{Lemmas}
Before proving the main theorems, we state and prove some lemmas that will be used in the proof.
In the first two lemmas, we treat $\trace \langle M \rangle_n$ and $\langle M \rangle_n$.
To this end, we collect properties of $\zeta_n$ and 
$
	\tau_n 
	= 
		\sum_{k=3}^n \zeta_k^{-2} 
$,
because \corref{col:Limit-of-conditional-expectation-and-covariance-of-Xn} allows us
to derive the asymptotic behavior of $\trace \langle M \rangle_n$ and $\langle M\rangle_n$
from that of $\zeta_n$ and $\tau_n$.
Note that \eqref{eq:Gamma-asymptotic} implies that $\zeta_n$ satisfies
\begin{align}\label{eq:zeta_n-asymptotic}
	\zeta_n
	=
		\left\{
			\prod_{k=2}^{n-1}
				\left(
					1+\frac{\delta}{k}
				\right)
		\right\}
		\zeta_2
	=
		\frac{\Gamma(n+\delta)}{\Gamma(n)}
		\frac{\Gamma(2)}{\Gamma(2+\delta)}
	\sim
		\frac{1}{\Gamma(2+\delta)}
		n^{\delta}
	\qquad
	\text{as}
	\qquad
	n\to\infty.
\end{align}
The behavior of $\tau_n$ is classified into three regimes according to the value of $\delta$. 
In each regime, by \eqref{eq:zeta_n-asymptotic}, $\tau_n$ satisfies the following asymptotics:
\begin{align}\label{eq:tau_n-order}
	&\tau_n
	\sim
		\begin{cases}
			\frac{\Gamma(2+\delta)^2}{1-2\delta}
			n^{1-2\delta} \quad &\text{if $\delta<\frac12$}, \\
			\frac{9\pi}{16}\log n \quad &\text{if $\delta=\frac12$}, \\
		\end{cases}
	&
	\lim_{n\to\infty} \tau_n <\infty \quad \text{if $\delta>\frac12$}.
\end{align}

We go back to $\trace \langle M \rangle_n$.
From \eqref{eq:Predictable-quad-variation-of-M_n}, we have
\begin{align}\label{eq:Trace-of-predictable-quadratic-variation}
	\trace \langle M \rangle_n
	=
		\sum_{k=3}^n \zeta_k^{-2} \trace (\covariance(X_k|\sigmaField_{k-1}))
	\leq
		C\tau_n.
\end{align}
Here $C$ is a positive constant.
The last estimate follows from
\begin{align*}
	\trace(\covariance(X_k|\sigmaField_{k-1}))
	\le
		\expect[|X_k|^2|\sigmaField_{k-1}]
	\le
		\sup_{v\in V} |v|^2.
\end{align*}
Then the following lemmas are shown.
\begin{lemma}\label{lem8423908409238104}
	The following holds.
	\begin{enumerate}
		\item	If $\delta \leq \frac{1}{2}$,
				then $\trace\langle M \rangle_n$ diverges almost surely as $n\to\infty$.
		\item	If $\delta > \frac{1}{2}$,
				then
				$
					\expect[\trace \langle M\rangle_\infty]
					<
						\infty
				$.
	\end{enumerate}
\end{lemma}
\begin{proof}
	The first assertion follows from 
	the former part of \eqref{eq:Trace-of-predictable-quadratic-variation},
	\eqref{eq:zeta_n-asymptotic},
	and \corref{col:Limit-of-conditional-expectation-and-covariance-of-Xn}.
	The second assertion follows from \eqref{eq:Trace-of-predictable-quadratic-variation}
	and \eqref{eq:tau_n-order}.
\end{proof}

\begin{lemma}\label{lem:<M>n/tau_n-Converges}
	If $\delta \le \frac{1}{2}$, then we have
	\begin{align*}
		\lim_{n\to\infty}
			\frac{\langle M\rangle_n}{\tau_n}
		=
			\frac{1}{m}
			\sum_{v\in V_\odd}
				(v - \bar{v})
				(v - \bar{v})^\top
		\quad
		\text{a.s.}
	\end{align*}
\end{lemma}
\begin{proof}
	From \eqref{eq:tau_n-order},
	$\zeta_k^{-2}/\tau_n$ converges to $0$ as $n\to\infty$ for each fixed $k$.
	Since $\sum_{k=3}^{n} \zeta_k^{-2}/\tau_n = 1$ for all $n$
	and the conditional covariance $\covariance(X_k|\sigmaField_{k-1})$ 
	converges almost surely to the constant matrix 
	from \corref{col:Limit-of-conditional-expectation-and-covariance-of-Xn},
	the Silverman-Toeplitz theorem implies the desired result.
\end{proof}

The next lemma gives an estimate of $R_n$.
\begin{lemma}\label{lem:Ln-estimation}
	For $n\geq 3$, we have
	\begin{align*}
		|R_n|
		\le
			\sum_{k=3}^n
				\frac{C}{k\zeta_k}
			+
			\frac{C}{\zeta_n}
		\leq
			C
			\begin{cases}
				1 & \text{if $\delta>0$}, \\
				\log n & \text{if $\delta=0$}, \\
				n^{-\delta} & \text{if $\delta<0$}.
			\end{cases}
	\end{align*}
	For $n \ge m\ge 3$, we have
	\begin{align*}
		|R_n - R_m|
		\le
			C
			\left(
				\frac{1}{\zeta_n}
				+\frac{1}{\zeta_m}
				+
					\sum_{k=m}^n \frac{1}{k\zeta_k}
			\right).
	\end{align*}
	Here, $C$ is a positive constant independent of $m$ and $n$.
\end{lemma}
\begin{proof}
	Write
	$
		r_k
		=
			R_k-R_{k-1}
	$
	for $k\geq 3$.
	Then $R_n  = \sum_{k=3}^n r_k$.
	Pairing the sum in $R_n$ by odd and even terms, we have
	\begin{align*}
		R_n
		=
			\sum_{l=2}^N
				(r_{2l-1} + r_{2l}) 
			+
			    \indicator{\odd}(n) r_n,
	\end{align*}
	where $N=\lfloor n/2 \rfloor$.
	Using \lref{lem:Conditional-expectation},
	we have
	$
		r_k
		=
			\zeta_k^{-1} ((-1)^{k-1}\frac{\gamma T_{k-1}}{k-1} + \xi_{k-1})
	$,
	where $\xi_i = \{(1-\gamma)(-1)^i - (\delta/i)\indicator{\odd}(i)\} \bar{v}$.
	Hence
	\begin{align*}
		r_{2l-1}+r_{2l}
		&=
			\zeta_{2l-1}^{-1} 
			\left(\frac{\gamma T_{2l-2}}{2l-2}+ \xi_{2l-2}\right)
			+
			\zeta_{2l}^{-1}
			\left(-\frac{\gamma T_{2l-1}}{2l-1}+ \xi_{2l-1}\right)\\
		&=
			\zeta_{2l}^{-1}
			\left(
				\frac{\gamma T_{2l-2}}{2l-2}
				+\xi_{2l-2}
				-\frac{\gamma T_{2l-1}}{2l-1}
				+\xi_{2l-1}
			\right)
			+
			\frac{\delta}{(2l-1)}
			\zeta_{2l}^{-1}
			\left(\frac{\gamma T_{2l-2}}{2l-2}+ \xi_{2l-2}\right).
	\end{align*}
	Here we used $\zeta_{2l-1}^{-1}=\zeta_{2l}^{-1}+\frac{\delta}{2l-1}\zeta_{2l}^{-1}$.
	Noting that $T_{2l-1} = T_{2l-2} + (-1)^{2l-2} X_{2l-1}$ 
	and $\xi_{2l-2} + \xi_{2l-1} = -\delta \bar{v}/(2l-1)$, this is further rewritten as
	\begin{align*}
		r_{2l-1}+r_{2l}
		=
			\frac{1}{(2l-1)\zeta_{2l}}
			\left(
				\frac{\gamma T_{2l-2}}{2l-2}
				-\gamma X_{2l-1}
				-\delta \bar{v}
			\right)
			+
			\frac{\delta}{(2l-1)\zeta_{2l}}
			\left(\frac{\gamma T_{2l-2}}{2l-2}+ \xi_{2l-2}\right).
	\end{align*}
	Hence using $|T_i| \le Ci$, 
	$|X_i| \le C$ and $|\xi_i| \le C$ for some constant $C>0$,
	we have
	\begin{align*}
		|r_{2l-1}+r_{2l}|
		\leq
			\frac{C}{(2l-1)\zeta_{2l}}
		\leq
			\frac{C}{(2l-1)\zeta_{2l-1}}
			+
			\frac{C}{2l\zeta_{2l}}.
	\end{align*}
	Using $|r_n|\leq C/\zeta_n$ and the above,
	\begin{align*}
		|R_n|
		\leq
			C
			\sum_{l=2}^N
				\left\{
					\frac{1}{(2l-1)\zeta_{2l-1}}
					+
					\frac{1}{2l\zeta_{2l}}
				\right\}
			+
			\frac{C}{\zeta_n}
		\leq
			\sum_{k=3}^n
				\frac{C}{k\zeta_k}
			+
			\frac{C}{\zeta_n}.
	\end{align*}
	The first assertion is proved.
	From the definition, we have
	\begin{align*}
		|R_n - R_m|
		=
			\left|
				\sum_{l=M+1}^N
						(r_{2l-1} + r_{2l}) 
					+
					\indicator{\odd}(n) r_n
					-
					\indicator{\odd}(m) r_m
			\right|
		\leq
			\sum_{l=M+1}^N
				|r_{2l-1} + r_{2l}|
			+
			|r_n|
			+
			|r_m|
			,
	\end{align*}
	where $M=\lfloor m/2 \rfloor$.
	Using the estimates above, we see the second assertion.
\end{proof}

\subsection{The diffusive regime}
We show the law of large numbers and the central limit theorem
in the diffusive regime ($\delta<1/2$).
\begin{proof}[Proof of \tref{thm:SLLN-in-diffusive}]
	Recall \lref{lem8423908409238104}, 
	\eqref{eq:Trace-of-predictable-quadratic-variation} and \eqref{eq:tau_n-order}.
	Then, it follows from \lref{lem:SLLN-for-trace-divergent-case} and \eqref{eq:tau_n-order}
	that for any $\eta>0$, 
	\begin{align*}
		|M_n|^2
		=
			o(\tau_n(\log \tau_n)^{1+\eta})
		=
			o(n^{1-2\delta}(\log n)^{1+\eta})
			 \quad \text{a.s.}
	\end{align*}
	Hence, since $S_n=\zeta_n\{M_n+R_n+S_2\}$, using the above, \lref{lem:Ln-estimation}
	and \eqref{eq:zeta_n-asymptotic}, we have
	\begin{align*}
		|\zeta_n M_n|
		&\leq
			C
			|n^\delta M_n|
		=
			o(\sqrt{n(\log n)^{1+\eta}}),
		&
		|\zeta_n R_n|
		&\leq
			C
			\begin{cases}
				n^\delta & \text{if $\delta>0$}, \\
				\log n & \text{if $\delta=0$}, \\
				1 & \text{if $\delta<0$}.
			\end{cases}
	\end{align*}
	Finally, using $|\zeta_n S_2|\leq Cn^\delta$ from \eqref{eq:zeta_n-asymptotic},
	we obtain the desired result.
\end{proof}

\begin{proof}[Proof of \tref{thm:CLT-in-diffusive}]
	We set $M_{n,k} = M_k/\sqrt{\tau_n}$ for $2\le k \le n$.
	Then, $\{M_{n,k}, \sigmaField_k\}_{2\le k \le n}$ is a mean-zero square-integrable martingale array with increments
	\begin{align*}
		\Delta M_{n,k} 
		=
			\frac{\Delta M_k}{\sqrt{\tau_n}}.
	\end{align*}
	To apply \lref{lem:CLT-for-vector-martingale-triangular-array},
	we will verify the convergence of the predictable quadratic variation and the Lindeberg condition.

	We compute the predictable quadratic variation of the martingale array $\{M_{n,k}\}$.
	From \eqref{eq:Predictable-quad-variation-of-M_n}, we have
	\begin{align*}
		\langle M_n \rangle_n
		=
			\sum_{k=3}^{n} \expect[\Delta M_{n,k} (\Delta M_{n,k})^\top |\sigmaField_{k-1}]
		=
			\frac{1}{\tau_n} \sum_{k=3}^{n} \expect[\Delta M_k (\Delta M_k)^\top |\sigmaField_{k-1}]
		=
			\frac{\langle M \rangle_n}{\tau_n}.
	\end{align*}
	From \lref{lem:<M>n/tau_n-Converges}, we obtain 
	\begin{align*}
		\lim_{n\to\infty}
			\langle M_n \rangle_n
		=
			\lim_{n\to\infty}
				\frac{\langle M\rangle_n}{\tau_n}
		=
			\frac{1}{m}
			\sum_{v\in V_\odd}
				(v - \bar{v})
				(v - \bar{v})^\top
		\quad
		\text{a.s.}
	\end{align*}

	Next, we verify the Lindeberg condition.
	Recalling \eqref{eq840928409231}, we have
	$|\Delta M_k|\leq C\zeta_k^{-1}$ for some $C>0$.
	Hence we obtain
	\begin{align*}
		|\Delta M_{n,k}|^2 \setindicator{|\Delta M_{n,k}| \ge \epsilon}
		=
			\frac{1}{\tau_n}
			|\Delta M_k|^2 \setindicator{|\Delta M_k| \ge \epsilon \sqrt{\tau_n}}
		\leq
			\frac{1}{(\epsilon \tau_n)^2}
			|\Delta M_k|^4
		\leq
			\frac{C}{(\epsilon \tau_n)^2}
			\frac{1}{\zeta_k^4},
	\end{align*}
	which implies
	\begin{align*}
		\sum_{k=3}^{n}
			\expect[|\Delta M_{n,k}|^2 \setindicator{|\Delta M_{n,k}| \ge \epsilon}|\sigmaField_{k-1}]
		&\leq
			\frac{C}{\epsilon^2}
			\frac{1}{\tau_n^2}
			\sum_{k=3}^{n} \frac{1}{\zeta_k^4}
		\leq
			\frac{C}{\epsilon^2}
			\begin{cases}
				n^{-1} & \text{if $\delta<1/4$}, \\
				n^{-1}\log n & \text{if $\delta=1/4$}, \\
				n^{-(2-4\delta)} & \text{if $1/4<\delta<1/2$}.
			\end{cases}
	\end{align*}
	In the last estimate, we used \eqref{eq:zeta_n-asymptotic} and \eqref{eq:tau_n-order}.
	For all cases of $\delta<1/2$, the right-hand side goes to $0$ as $n\to\infty$.
	Hence, the Lindeberg condition holds.
	Therefore, we conclude from \lref{lem:CLT-for-vector-martingale-triangular-array} that
	\begin{align*}
		\lim_{n\to\infty}
			M_{n,n}
		&=
			\mathcal{N}_d
			\left(\,\zerovector_d, \; \frac{1}{m}\sum_{v\in V_\odd}
				(v - \bar{v})
				(v - \bar{v})^\top\right)
		\qquad
		\text{in distribution.}
	\end{align*}
	Hence, recalling $\frac{S_n}{\sqrt{n}}=\frac{\zeta_n}{\sqrt{n}}\{\sqrt{\tau_n}M_{n,n}+R_n+S_2\}$
	and using
	\begin{align*}
		\lim_{n\to\infty}
			\frac{\zeta_n\sqrt{\tau_n}}{\sqrt{n}}
		&=
			\frac{1}{\sqrt{1-2\delta}},
		&
		\lim_{n\to\infty}
			\frac{\zeta_n R_n}{\sqrt{n}}
		&=
			0,
		&
		\lim_{n\to\infty}
			\frac{\zeta_n S_2}{\sqrt{n}}
		&=
			0
		\qquad
		\text{a.s.},
	\end{align*}
	which follows from \eqref{eq:zeta_n-asymptotic}, \eqref{eq:tau_n-order} and \lref{lem:Ln-estimation},
	we arrive at the assertion.
\end{proof}

\subsection{The critical regime}

Next, we prove the law of large numbers and the central limit theorem
in the critical regime ($\delta=1/2$).
\begin{proof}[Proof of \tref{thm:SLLN-in-critical}]
	Recall \lref{lem8423908409238104}, 
	\eqref{eq:Trace-of-predictable-quadratic-variation} and \eqref{eq:tau_n-order}.
	Then, from \lref{lem:SLLN-for-trace-divergent-case}, 
	for any $\eta>0$, we have $|M_n|=o(\sqrt{(\log n)(\log \log n)^{1+\eta}})$ almost surely.
	Recall $S_n=\zeta_n\{M_n+R_n+S_2\}$. 
	Noting that $|R_n|$ and $|S_2|$ are bounded from \lref{lem:Ln-estimation}
	and using the asymptotic behavior of $\zeta_n$ in \eqref{eq:zeta_n-asymptotic},
	we obtain the desired result.
\end{proof}

\begin{proof}[Proof of \tref{thm:CLT-in-critical}]
	As in the proof of the diffusive regime, we apply \lref{lem:CLT-for-vector-martingale-triangular-array} 
	to the martingale array $\{M_{n,k}\}$ defined in the proof of \tref{thm:CLT-in-diffusive}.
	The convergence of the predictable quadratic variation is obtained from \lref{lem:<M>n/tau_n-Converges} as before.
	To verify the Lindeberg condition, 
	we observe  from \eqref{eq:zeta_n-asymptotic} and \eqref{eq:tau_n-order} that
	\begin{align*}
		\frac{1}{\tau_n^2}
		\sum_{k=3}^{n}
			\frac{1}{\zeta_k^4}
		\leq
			\frac{C}{(\log n)^2}
			\sum_{k=3}^{n}
				\frac{1}{k^2}
		\to
			0.
	\end{align*}
	By the same argument as in the proof of the diffusive regime, the Lindeberg condition holds.
	Therefore, we conclude from \lref{lem:CLT-for-vector-martingale-triangular-array} that
	\begin{align*}
		\lim_{n\to\infty}
			M_{n,n}
		=
			\mathcal{N}_d\left(\,\zerovector_d, \; 
				\frac{1}{m}
				\sum_{v\in V_\odd}
					(v - \bar{v})
					(v - \bar{v})^\top\right)
		\qquad
		\text{in distribution}.
	\end{align*}
	From \eqref{eq:zeta_n-asymptotic}, \eqref{eq:tau_n-order} and \lref{lem:Ln-estimation}, we have
	\begin{align*}
		\lim_{n\to\infty}
			\frac{\zeta_n \sqrt{\tau_n}}{\sqrt{n\log n}}
		&=
			1,
		&
		\lim_{n\to\infty}
			\frac{\zeta_n R_n}{\sqrt{n\log n}}
		&=
			0,
		&
		\lim_{n\to\infty}
			\frac{\zeta_n S_2}{\sqrt{n\log n}}
		&=
			0
		\quad \text{a.s.}
	\end{align*}
	Noting $\frac{S_n}{\sqrt{n\log n}}=\frac{\zeta_n}{\sqrt{n\log n}}\{\sqrt{\tau_n}M_{n,n}+R_n+S_2\}$,
	we arrive at the assertion.
\end{proof}

\subsection{The superdiffusive regime}

Finally, we prove the law of large numbers in the superdiffusive regime ($\delta>1/2$).
\begin{proof}[Proof of \tref{thm:SLLN-in-superdiffusive}]
	\lref[lem8423908409238104]{lem:SLLN-for-trace-convergent-case} imply that
	$
		M_n \to M_\infty 
	$
	almost surely and in $L^2$ as $n\to\infty$.
	Since $M_n=\zeta_n^{-1} S_n - R_n - S_2$, we show that $R_n$ converges almost surely and in $L^2$.
	From \eqref{eq:Gamma-asymptotic}, for $\delta>0$, we have
	\begin{align*}	
		\lim_{m,n\to\infty}
			\sum_{k=m}^n \frac{1}{k\zeta_k}
		=
			0.
	\end{align*}
	Since $\delta>1/2$, this summation and $\zeta_n^{-1}, \zeta_{m}^{-1}$ go to $0$ as $n,m \to \infty$.
	Thus, with the help of \lref{lem:Ln-estimation},
	$\{R_n\}_{n\geq 2}$ is a Cauchy sequence almost surely.
	Hence $R_n$ converges to $R_\infty$ almost surely.
	Since $\{R_n\}$ is uniformly bounded, the bounded convergence theorem
	implies that $\{R_n\}$ converges in $L^2$ as well.
	Noting that $\zeta_n^{-1} S_n = M_n + R_n + S_2$, these results and \eqref{eq:zeta_n-asymptotic} yield that
	\begin{align*}
		\lim_{n\to\infty}
			\frac{S_n}{n^\delta}
		=
			\lim_{n\to\infty}
				\frac{\zeta_n}{n^\delta}
				\frac{S_n}{\zeta_n}
		=
			\frac{M_\infty+R_\infty+S_2}{\Gamma(2+\delta)}.
	\end{align*}
	Thus, the first identity is shown.
	From \tref{thm43092043241902494}, 
	\begin{align*}
		\expect[L]
		=
			\lim_{n\to\infty}
				n^{-\delta}
				\expect[S_n]
		=
			\zerovector_d.
	\end{align*}
	We also have
	\begin{align*}
		\expect[|L|^2]
		=
			\expect[|L - \expect[L]|^2]
		=
			\lim_{n\to\infty}
				n^{-2\delta}
				\expect
					[
						|S_n-\expect[S_n]|^2
					]
		=
			C_{\alpha,\beta}.
	\end{align*}
	Here we used \tref{thm:Asymptotic-second-moment}.
	The proof is complete.
\end{proof}

\section{Concluding remarks}\label{secConcludingRemarks}
We conclude this paper by making comments on the cases where $V_\odd\cap V_\even\not=\emptyset$.
Under \aref{assumption49023492104} without Assumption \eqref{item500},
we can define the ERW. Indeed, we modify the distribution of $G^\odd_k(v)$ as follows.
Let $v \in V$ and $w \in V_\odd$.
\begin{itemize}
	\item If $v \in V_\odd$, then we define
		\begin{align*}
			\prob(G^\odd_k(v) = w)
			=
				\begin{cases}
					p & \text{if $w = v$}, \\
					\frac{1-p}{m-1} & \text{if $w \neq v$}.
				\end{cases}
		\end{align*}
	\item If $v \in V \setminus V_\odd$ and hence $- v \in V_\odd$, then we define
		\begin{align*}
			\prob(G^\odd_k(v) = w)
			=
				\begin{cases}
					q & \text{if $ w = -v $}, \\
					\frac{1-q}{m-1} & \text{if $ w \neq -v $}.
				\end{cases}
		\end{align*}
		We do not consider $\prob(G^\odd_k(v) = w)$ if $V\setminus V_\odd$ is empty. 
\end{itemize}
Similarly, the definition of $G^\even_k(v)$ is also changed.
Then the ERW is well-defined.

First we consider the case $V_\odd=V_\even$.
As in the proof of \lref{lem:Conditional-expectation}, we have
$
	A^+_n
	=
		2\alpha S_n
$
and
$
	A^-_n
	=
		2(1-\alpha)n\bar{v}
$.
Hence, since $\bar{v}=0$ holds from $V_\odd=-V_\even$ and $V_\odd=V_\even$, we have
\begin{align*}
	\expect[X_{n+1} | \sigmaField_n]
	=
		\frac{\alpha}{n} S_n.
\end{align*}
This coincides with the formula for the conditional expectation of the step of the ERW studied in \cite{BercuLaulin2019MERW},
so the identical case falls within the framework of \cite{BercuLaulin2019MERW}.
In this framework, we can treat 
\exref{ex84294820942} \itemref{item890284309194}
and the usual triangular lattice.

In the partially overlapping case, that is, 
when $V_\odd \cap V_\even \neq \emptyset$ and $V_\odd \neq V_\even$,
it is difficult to analyze the ERW for the reason stated in the introduction.
Although an expression for $\expect[X_{n+1} | \sigmaField_n]$ can still be derived,
we do not obtain an analogue of \tref{thm:Asymptotic-second-moment}, 
and so the regimes cannot be defined in the same way.
Note that the directions in $V_\odd\cap V_\even$
contribute at every step, whereas the directions in 
$V_\odd\setminus(V_\odd\cap V_\even)=V_\odd\setminus V_\even$
and $V_\even\setminus(V_\odd\cap V_\even)=V_\even\setminus V_\odd$
appear only when the walk is located in the corresponding local structure.
Because of this asymmetry, we do not know how to analyze this case with our present methods.
Therefore, the asymptotic behavior in this case remains open.

\appendix

\section{Limit theorems for vector martingales}\label{appendix}
In this section, we provide strong laws of large numbers and a central limit theorem for vector martingales.
The results presented in this section are standard and can be found, for instance, 
in Duflo \cite{Duflo1997Book}.

Let $\{M_n\}_{n\ge0}$ be a square integrable $\RealNum^d$-valued martingale with $M_0=0$. 
Write $\Delta M_n = M_n - M_{n-1}$ for $n\ge1$.
The predictable quadratic variation associated with $M_n$ is given, for all $n\ge1$, by
\begin{align}\label{eq49320941032413}
	\langle M \rangle_n
	=
		\sum_{k=1}^{n} \expect\left[\Delta M_k (\Delta M_k)^\top|\sigmaField_{k-1}\right].
\end{align}
It is easy to see that $\{M_n M_n^\top - \langle M \rangle_n\}_{n\ge0}$ is a martingale,
namely, it holds that
\begin{align}\label{eq4829084209}
	\expect[M_{n+1} M_{n+1}^\top - \langle M \rangle_{n+1}|\sigmaField_n]
	=
		M_n M_n^\top - \langle M \rangle_n.
\end{align}
In particular, we have
$
	\expect[M_n M_n^\top]=\expect[\langle M \rangle_n]
$
and 
$
	\expect[|M_n|^2]=\expect[\trace \langle M \rangle_n]
$.

The following lemmas provide an asymptotic estimate and limit theorem of $M_n$.
\begin{lemma}\label{lem:SLLN-for-trace-divergent-case}
	For any $\eta>0$, we have
	\begin{align*}
		|M_n|^2 = o(\trace \langle M \rangle_n (\log \trace \langle M \rangle_n)^{1+\eta})
		\quad
		\text{a.s. on\, $\{\trace \langle M \rangle_\infty = \infty\}$.}
	\end{align*}
\end{lemma}
\begin{lemma}\label{lem:SLLN-for-trace-convergent-case}
	If $\expect[\trace \langle M \rangle_\infty] < \infty$, then we have
	\begin{align*}
		\lim_{n\to\infty}
			M_n
		=
			M_\infty
		\quad \text{a.s. and\, in $L^2$.}
	\end{align*}
\end{lemma}

Next, we give a central limit theorem for vector martingales.

\begin{lemma}\label{lem:CLT-for-vector-martingale-triangular-array}
	Let $\{k_n\}_{n\ge1}$ be a sequence of positive integers with $k_n\to\infty$.
	For each $n\ge1$, let $\{\sigmaField_{n,i}\}_{0\le i\le k_n}$ be a filtration and let
	$\{M_{n,i}\}_{0\le i\le k_n}$ be a square integrable $\RealNum^d$-valued martingale with $M_{n,0}=0$
	adapted to $\{\sigmaField_{n,i}\}$.
	Assume the following.
	\begin{enumerate}
		\item	There exists a deterministic symmetric positive semidefinite matrix
				$\Sigma\in \RealNum^{d\times d}$ such that
				$
					\langle M_n\rangle_{k_n}
					\to
						\Sigma
				$
				in probability.
				Here $\langle M_n\rangle$ is
				the predictable quadratic variation associated with $M_n$
				defined by \eqref{eq49320941032413}.
		\item	For any $\epsilon>0$,
				$
					\sum_{i=1}^{k_n}
						\expect\left[
							|\Delta M_{n,i}|^2 \setindicator{|\Delta M_{n,i}|>\epsilon}
							\mid \sigmaField_{n,i-1}
						\right]
					\to
						0
				$
				in probability.
	\end{enumerate}
	Then, $\{M_{n,k_n}\}_{n\geq 1}$ converges to $\mathcal{N}_d(0,\Sigma)$
	in distribution.
\end{lemma}


\begin{thebibliography}{10}

\bibitem{BaurBertoin2016ERWUrn}
E.~Baur and J.~Bertoin.
\newblock Elephant random walks and their connection to {P}\'olya-type urns.
\newblock {\em Phys. Rev. E}, 94:052134, 2016.

\bibitem{Bercu2018MartingaleERW}
B.~Bercu.
\newblock A martingale approach for the elephant random walk.
\newblock {\em J. Phys. A}, 51(1):015201, 16, 2018.

\bibitem{Bercu2025MERWWithStops}
B.~Bercu.
\newblock On the multidimensional elephant random walk with stops.
\newblock {\em Stochastic Process. Appl.}, 189:Paper No. 104692, 17, 2025.

\bibitem{BercuLaulin2019MERW}
B.~Bercu and L.~Laulin.
\newblock On the multi-dimensional elephant random walk.
\newblock {\em J. Stat. Phys.}, 175(6):1146--1163, 2019.

\bibitem{Bertenghi2022FunctionalMERW}
M.~Bertenghi.
\newblock Functional limit theorems for the multi-dimensional elephant random walk.
\newblock {\em Stoch. Models}, 38(1):37--50, 2022.

\bibitem{ChenLaulin2023MARW}
J.~Chen and L.~Laulin.
\newblock Analysis of the smoothly amnesia-reinforced multidimensional elephant random walk.
\newblock {\em J. Stat. Phys.}, 190(10):Paper No. 158, 42, 2023.

\bibitem{ColettiGavaSchutz2017CLT}
C.~F. Coletti, R.~Gava, and G.~M. Sch\"utz.
\newblock Central limit theorem and related results for the elephant random walk.
\newblock {\em J. Math. Phys.}, 58(5):053303, 8, 2017.

\bibitem{ColettiGavaSchutz2017SIP}
C.~F. Coletti, R.~Gava, and G.~M. Sch\"utz.
\newblock A strong invariance principle for the elephant random walk.
\newblock {\em J. Stat. Mech. Theory Exp.}, (12):123207, 8, 2017.

\bibitem{CurienLaulin2024PlaneRecurrenceMERW}
N.~Curien and L.~Laulin.
\newblock Recurrence of the plane elephant random walk.
\newblock {\em C. R. Math. Acad. Sci. Paris}, 362:1183--1188, 2024.

\bibitem{Duflo1997Book}
M.~Duflo.
\newblock {\em Random iterative models}, volume~34 of {\em Applications of Mathematics (New York)}.
\newblock Springer-Verlag, Berlin, 1997.
\newblock Translated from the 1990 French original by Stephen S. Wilson and revised by the author.

\bibitem{GhoshDhillonKataria2026MovesMERWS}
S.~Ghosh, M.~Dhillon, and K.~K. Kataria.
\newblock On limiting behaviour of moves in multidimensional elephant random walk with stops.
\newblock {\em arXiv:2601.07502}, 2026.

\bibitem{GonzalezNavarrete2020RandomTendencyMERW}
M.~Gonz\'alez-Navarrete.
\newblock Multidimensional walks with random tendency.
\newblock {\em J. Stat. Phys.}, 181(4):1138--1148, 2020.

\bibitem{GuerinLaulinRaschel2023FixedPointSuperdiffusiveMERW}
H.~Gu\'erin, L.~Laulin, and K.~Raschel.
\newblock A fixed-point equation approach for the superdiffusive elephant random walk.
\newblock {\em Ann. Inst. Henri Poincar{\'e}, Probab. Stat.}, 62(2): 973--1005, 2026.

\bibitem{KotaniSunada2000Jacobian}
M.~Kotani and T.~Sunada.
\newblock Jacobian tori associated with a finite graph and its abelian covering graphs.
\newblock {\em Adv. in Appl. Math.}, 24(2):89--110, 2000.

\bibitem{KotaniSunada2001Standard}
M.~Kotani and T.~Sunada.
\newblock Standard realizations of crystal lattices via harmonic maps.
\newblock {\em Trans. Amer. Math. Soc.}, 353(1):1--20, 2001.

\bibitem{KotaniSunada2003}
M.~Kotani and T.~Sunada.
\newblock Spectral geometry of crystal lattices.
\newblock In {\em Heat kernels and analysis on manifolds, graphs, and metric spaces ({P}aris, 2002)}, volume 338 of {\em Contemp. Math.}, pages 271--305. Amer. Math. Soc., Providence, RI, 2003.

\bibitem{Marquioni2019CoupledMemoryMERW}
V.~M. Marquioni.
\newblock Multidimensional elephant random walk with coupled memory.
\newblock {\em Phys. Rev. E}, 100(5):052131, 10, 2019.

\bibitem{Mukherjee2025CayleyTree}
S.~S. Mukherjee.
\newblock Elephant random walks on infinite Cayley trees.
\newblock {\em arXiv:2509.03048}, 2025.

\bibitem{Qin2025RecTransMERW}
S.~Qin.
\newblock Recurrence and transience of multidimensional elephant random walks.
\newblock {\em Ann. Probab.}, 53(3):1049--1078, 2025.

\bibitem{SchutzTrimper2004ERW}
G.~M. Sch\"utz and S.~Trimper.
\newblock Elephants can always remember: Exact long-range memory effects in a non-{M}arkovian random walk.
\newblock {\em Phys. Rev. E}, 70:045101, 2004.

\bibitem{Shibata2025ArXivPeriodic}
S.~Shibata.
\newblock Functional limit theorems for elephant random walks on general periodic structures.
\newblock {\em arXiv:2511.10347}, 2025.

\bibitem{Sunada2012}
T.~Sunada.
\newblock Lecture on topological crystallography.
\newblock {\em Jpn. J. Math.}, 7(1):1--39, 2012.

\bibitem{Sunada2013Book}
T.~Sunada.
\newblock {\em Topological crystallography. With a view towards discrete geometric analysis}, volume~6 of {\em Surveys and Tutorials in the Applied Mathematical Sciences}.
\newblock Springer, Tokyo, 2013.

\bibitem{Ureshino2021}
Y.~Ureshino.
\newblock Elephant random walk on multi-dimensional lattices.
\newblock Master's thesis, Keio University, 2021 (in Japanese).

\bibitem{YanChengBai2006}
X.~Yan, Y.~Cheng, and Z.~Bai.
\newblock Asymptotics of adaptive design with two alternating generating matrices.
\newblock {\em J. Statist. Plann. Inference}, 136(11):4043--4058, 2006.

\end{thebibliography}

\end{document}